\documentclass{article}

\usepackage[titletoc,title]{appendix}
\usepackage[cp1250]{inputenc}
\usepackage[IL2]{fontenc}
\usepackage{yfonts,fancyhdr}
\usepackage{a4wide}
\usepackage[english]{babel}
\usepackage{euscript}
\usepackage{amstext,amsbsy,amscd,amssymb}
\usepackage{amsmath,enumerate}
\usepackage{amsfonts}
\usepackage{mathrsfs}
\usepackage{graphics}
\usepackage{microtype}
\allowdisplaybreaks

\let\rarr=\rightarrow

\let\veps=\varepsilon
\let\mcal=\mathcal
\let\mfrak=\mathfrak
\let\eus=\EuScript

\def\N{\mathbb{N}}
\def\Z{\mathbb{Z}}
\def\R{\mathbb{R}}
\def\C{\mathbb{C}}

\def\Hom{\mathop {\rm Hom} \nolimits}

\def\Ind{\mathop {\rm Ind} \nolimits}

\def\ad{\mathop {\rm ad} \nolimits}

\def\Ad{\mathop {\rm Ad} \nolimits}

\def\Sol{\mathop {\rm Sol} \nolimits}

\def\tr{\mathop {\rm tr} \nolimits}

\def\ch{\mathop {\rm ch} \nolimits}

\newcommand{\mC}{\mathbb C}

\newcommand{\mN}{\mathbb N}

\newsymbol\squares 1003

\long\def\proof #1{\noindent \emph{Proof.}\ #1 \hfill $\squares$
\medskip}

\newcounter{num}[section]
\numberwithin{equation}{section}
\numberwithin{num}{section}

\long\def\definition #1 {\refstepcounter{num} \noindent {\bf
Definition \thenum.} #1

\medskip}

\long\def\theorem #1{\refstepcounter{num} \noindent {\bf Theorem
\thenum.} #1

\medskip}

\long\def\lemma #1{\refstepcounter{num}  \noindent {\bf Lemma
\thenum.} #1

\medskip}

\long\def\proposition #1{\refstepcounter{num}  \noindent {\bf
Proposition \thenum.} #1

\medskip}

\long\def\remark{\noindent {\bf Remark.}\ }

\newenvironment{enum}{\begin{list}{}{\topsep=2pt \itemsep=0pt
\parsep=0pt}}{\end{list}}

\newcommand*\riso{%
  \xrightarrow[]{\raisebox{-0.25em}{\smash{\ensuremath{\sim}}}}%
}

\makeatletter
\newcommand*\if@single[3]{%
  \setbox0\hbox{${\mathaccent"0362{#1}}^H$}%
  \setbox2\hbox{${\mathaccent"0362{\kern0pt#1}}^H$}%
  \ifdim\ht0=\ht2 #3\else #2\fi
  }
\newcommand*\rel@kern[1]{\kern#1\dimexpr\macc@kerna}
\newcommand*\widebar[1]{\@ifnextchar^{{\wide@bar{#1}{0}}}{\wide@bar{#1}{1}}}
\newcommand*\wide@bar[2]{\if@single{#1}{\wide@bar@{#1}{#2}{1}}{\wide@bar@{#1}{#2}{2}}}
\newcommand*\wide@bar@[3]{%
  \begingroup
  \def\mathaccent##1##2{%
    \if#32 \let\macc@nucleus\first@char \fi
    \setbox\z@\hbox{$\macc@style{\macc@nucleus}_{}$}%
    \setbox\tw@\hbox{$\macc@style{\macc@nucleus}{}_{}$}%
    \dimen@\wd\tw@
    \advance\dimen@-\wd\z@
    \divide\dimen@ 3
    \@tempdima\wd\tw@
    \advance\@tempdima-\scriptspace
    \divide\@tempdima 10
    \advance\dimen@-\@tempdima
    \ifdim\dimen@>\z@ \dimen@0pt\fi
    \rel@kern{0.6}\kern-\dimen@
    \if#31
      \overline{\rel@kern{-0.6}\kern\dimen@\macc@nucleus\rel@kern{0.4}\kern\dimen@}%
      \advance\dimen@0.4\dimexpr\macc@kerna
      \let\final@kern#2%
      \ifdim\dimen@<\z@ \let\final@kern1\fi
      \if\final@kern1 \kern-\dimen@\fi
    \else
      \overline{\rel@kern{-0.6}\kern\dimen@#1}%
    \fi
  }%
  \macc@depth\@ne
  \let\math@bgroup\@empty \let\math@egroup\macc@set@skewchar
  \mathsurround\z@ \frozen@everymath{\mathgroup\macc@group\relax}%
  \macc@set@skewchar\relax
  \let\mathaccentV\macc@nested@a
  \if#31
    \macc@nested@a\relax111{#1}%
  \else
    \def\gobble@till@marker##1\endmarker{}%
    \futurelet\first@char\gobble@till@marker#1\endmarker
    \ifcat\noexpand\first@char A\else
      \def\first@char{}%
    \fi
    \macc@nested@a\relax111{\first@char}%
  \fi
  \endgroup
}
\makeatother

\newcommand\rsmraise[1]{%
  \ifx#1\displaystyle .8\else
    \ifx#1\textstyle .8\else
      \ifx#1\scriptstyle .6\else
        .45%
      \fi
    \fi
  \fi}

\title{Conformal Galilei algebras, symmetric polynomials and singular vectors}

\author{Libor Křižka, Petr Somberg}

\AtEndDocument{\bigskip{\footnotesize%
  (L.\,Křižka) \textsc{Charles University in Prague, Faculty of Mathematics and Physics, Mathematical Institute of Charles University, Sokolovská 83, 180\,00 Praha 8, Czech Republic} \par
  \textit{E-mail address}: \texttt{krizka@karlin.mff.cuni.cz} \par
  \addvspace{\medskipamount}
  (P.\,Somberg) \textsc{Charles University in Prague, Faculty of Mathematics and Physics, Mathematical Institute of Charles University, Sokolovská 83, 180\,00 Praha 8, Czech Republic} \par
  \textit{E-mail address}: \texttt{somberg@karlin.mff.cuni.cz} \par
}}

\begin{document}
\date{}
\maketitle

\begin{abstract}

We classify and explicitly describe homomorphisms of Verma modules for conformal Galilei algebras
$\mathfrak{cga}_\ell(d,\C)$ with $d=1$ for any integer value $\ell \in \mathbb{N}$. The homomorphisms
are uniquely determined by singular vectors as solutions of certain differential operators of flag
type, and identified with specific polynomials arising as coefficients in the expansion of a parametric
family of symmetric polynomials into power sum symmetric polynomials.


\end{abstract}

\thispagestyle{empty}

\tableofcontents


\section*{Introduction}
\addcontentsline{toc}{section}{Introduction}

An important structure associated to a system of partial
differential equations on manifolds is its algebra of symmetries.
In the case when such a system is of geometric origin, this
algebra allows to extract many fundamental properties of the
former system (e.g.\ generates from a given solution of the system
other solutions) and has a neat and delicate relationship to various properties
of the underlying manifold.

As for the perspective of their applications in theoretical physics,
one of the most important examples of such a system of partial differential
equations is the free Schrödinger equation. A classical result says that its symmetry
algebra is recognized as the Schrödinger algebra, which is a basic
representative of the family of the so called conformal Galilei algebras, see \cite{Havas-Plebanski1978}. The real conformal Galilei algebras $\mfrak{cga}_\ell(d,\R)$ parameterized by $d \in \N$ (called the spatial dimension) and $\ell \in {1 \over 2}\N_0$
(called the spin), are (generalizations of) non-relativistic analogues of Lie algebras
$\mathfrak{so}(d+1,2,\R)$ acting by conformal transformations on $\R^{d,1}$
and constructed as their contractions. In contrast with conformal Lie algebras, the conformal
Galilei algebras are not semisimple. There are several topics of recent interest related to conformal Galilei algebras, e.g.\ the classical mechanics with higher order time derivatives, the non-relativistic AdS/CFT correspondence for conformal Galilei algebras related to the flat holography in the limit when the radius of the AdS spacetime tends to infinity.

A standard mathematical setup consists in the change of the previous perspective. Having the symmetry realized by a conformal Galilei algebra, we may ask for the classification of its invariants. In the present article, we consider the pair consisting of a conformal Galilei algebra $\mathfrak{cga}_\ell(d,\C)$ for $d=1$ with any integer spin $\ell \in \N_0$ and its standard Borel subalgebra, and analyze the composition structure of Verma modules induced from characters of the standard Borel subalgebra. Homomorphisms of Verma modules are determined by the highest/lowest weight (or singular) vectors, which have another important significance. Namely, they are in a bijection with invariant differential operators acting on smooth sections of principal series representations supported on the homogeneous space for the conformal Galilei group and
induced from the dual character to that used for Verma modules. Although our considerations on the homogeneous space for conformal Galilei group are purely local, on a general smooth manifold plays the role of the transitive action of conformal Galilei group the Newton-Cartan geometric structure.

Based on the techniques of geometric representation theory and harmonic analysis, the main achievement in the present article is the complete classification and precise positions of singular vectors in Verma modules for conformal Galilei algebras $\mathfrak{cga}_\ell(1,\C)$ with $\ell \in \N_0$, thereby completing several partial and non-complete results in this direction we are aware of, see \cite{Aizawa-Dobrev-Doebner2002}, \cite{Aizawa-Isaac-Kimura2012}, \cite{Aizawa-Kimura-Segar2013}, \cite{Aizawa-Chandrashekar-Segar2015}. The complete classification of singular vectors in Verma modules for conformal Galilei algebras $\mathfrak{cga}_\ell(1,\C)$ with $\ell \in {1 \over 2} + \N_0$ is given in \cite{Aizawa-Kimura-Segar2013}. The starting point leading to our results is based on the reformulation of the former task to a problem of classifying polynomial solutions of certain differential equations of flag type. We found it quite convenient and rather indispensable to express all solutions of this equation in terms of certain polynomials organized into finite products of generating formal power series.
Let us note that a classification of all finite weight modules over the conformal Galilei algebras $\mfrak{cga}_\ell(1,\C)$ with any $\ell \in {1 \over 2}\N_0$ has been done in \cite{Lu-Mazorchuk-Zhao2014}.

Let us briefly summarize the content of our article. In Section \ref{sec:geomrealindmod},
we describe the geometric realization of induced modules on the homogeneous space for a general complex finite-dimensional Lie group and its Borel subgroup. In Section \ref{sec:confgalalg}, we discuss representation theoretical properties of conformal Galilei algebras in the spatial dimension $d=1$ together with a class of lowest weight modules. Namely, we describe $\mcal{D}$-module realization of these modules supported on the closed (point) orbit for the Borel subgroup. In fact, they are realized as the dual spaces of jets of sections of line bundles supported at the unit coset of the homogeneous space for conformal Galilei group. The complete classification together with explicit formulas for singular vectors is achieved by passing to the
Fourier dual picture, see Section \ref{sec:vermmodsingsec}. This converts the
former problem to describe the solution space of the differential
equation of flag type. The rest of this section is then devoted to an explicit description
of its solution space through the technique of generating formal power series and their
products. This is quite convenient and also surprising way due to a canonical
multiplicative structure on the solution space. For the reader's convenience,
we remind in Appendix \ref{app:symetric polynomials and fucntions} a few
basic formulas related to symmetric polynomials which are needed in the discussion.
We close the article by highlighting several open questions and problems related to our
work.

Throughout the article, we use the standard notation $\mathbb{Z}$, $\mathbb{N}$ and
$\mathbb{N}_0$ for the set of integers, the set of natural numbers and the set
of natural numbers including zero, respectively.


\section{Geometric realization of induced modules}
\label{sec:geomrealindmod}

In this section we describe a geometric realization of Verma modules for finite-dimensional complex Lie algebras, which are not necessarily complex semisimple Lie algebras.

Let $G$ be a complex finite-dimensional Lie group with its Lie algebra $\mfrak{g}$. Furthermore, let $\mfrak{h}$ be a Cartan subalgebra of $\mfrak{g}$, which means that $\mfrak{h}$ is a nilpotent Lie subalgebra of $\mfrak{g}$ such that
$\mfrak{n}_\mfrak{g}(\mfrak{h})=\mfrak{h}$, where $\mfrak{n}_\mfrak{g}(\mfrak{h})$ denotes the normalizer of $\mfrak{h}$ in $\mfrak{g}$. Let us note that there exists a Cartan subalgebra of $\mfrak{g}$ and that two Cartan subalgebras of $\mfrak{g}$ are conjugate by an automorphism of $\mfrak{g}$, and in particular they are isomorphic.
If we define a generalized eigenspace $\mfrak{g}_\alpha \subset \mfrak{g}$ for $\alpha \in \mfrak{h}^*$ by
\begin{align}
  \mfrak{g}_\alpha = \{X \in \mfrak{g};\, (\ad(H)-\alpha(H))^kX=0\ \text{for all $H \in \mfrak{h}$ and $k\gg 0$}\},
\end{align}
then we may introduce the notion of root spaces, roots, positive and negative roots for the Lie algebra $\mfrak{g}$ as
in the case of semisimple Lie algebras.
We denote by $\Delta$ the root system of $\mfrak{g}$ with respect to $\mfrak{h}$ and by $\Delta^+$ and $\Delta^-$ a positive and negative root subspaces in $\Delta$, respectively.

We associate to the positive and negative root system $\Delta^+$ and $\Delta^-$ the nilpotent Lie subalgebras
\begin{align}
  \mfrak{n} = \bigoplus_{\alpha \in \Delta^+} \mfrak{g}_\alpha \qquad \text{and} \qquad \widebar{\mfrak{n}} = \bigoplus_{\alpha \in \Delta^-} \mfrak{g}_\alpha,
\end{align}
and solvable Lie subalgebras
\begin{align}
  \mfrak{b}= \mfrak{h} \oplus \mfrak{n} \qquad \text{and} \qquad \widebar{\mfrak{b}}= \mfrak{h} \oplus \widebar{\mfrak{n}}
\end{align}
of $\mfrak{g}$, respectively. Moreover, we get a triangular decomposition
\begin{align}
  \mfrak{g} = \widebar{\mfrak{n}} \oplus \mfrak{h} \oplus \mfrak{n} \label{eq:triangular decomposition ex}
\end{align}
of the Lie algebra $\mfrak{g}$.
We refer the reader to \cite{Jacobson1979} for more detailed exposition of Cartan subalgebras of general Lie algebras.
Let us note that the Lie subalgebras $\mfrak{b}$ and $\widebar{\mfrak{b}}$ are not in general Borel subalgebras in the sense of maximal solvable Lie subalgebras of $\mfrak{g}$.
However, we will call $\mfrak{b}$ and $\widebar{\mfrak{b}}$ the standard and opposite standard Borel subalgebra of $\mfrak{g}$, respectively, because of the analogy with the triangular decomposition of semisimple Lie algebras.

Since $\mfrak{n}_\mfrak{g}(\mfrak{b})=\mfrak{b}$, we define the Borel subgroup $B$ of $G$ with its Lie algebra $\mfrak{b}$ as the connected component of the identity of the closed subgroup $N_G(\mfrak{b})$ of $G$. Then $p \colon G \rarr G/B$ is a principal $B$-bundle and $X=G/B$ is the corresponding flag manifold for $G$. We define the closed subgroup $\widebar{N}$ of $G$ to be the image of $\widebar{\mfrak{n}}$ under the exponential mapping $\exp \colon \mfrak{g} \rarr G$, so that the mapping
\begin{align}
  \exp \colon \widebar{\mfrak{n}} \rarr \widebar{N}
\end{align}
is a diffeomorphism. Therefore, for $g \in G$ the mapping $f_g\colon \widebar{N} \rarr X$ given by
\begin{align}
  f_g(n)=p(gn)=gnB
\end{align}
for $n \in \widebar{N}$ is an open embedding, and the manifold $X$ is covered by open subsets $U_g=p(g\widebar{N})$
for which we have
\begin{align}
  X= \bigcup_{g \in G} U_g.
\end{align}
Hence, we obtain the atlas $\{(U_g,u_g)\}_{g \in G}$ on $X$, where $u_g \colon U_g \rarr \widebar{\mfrak{n}}$ is defined by
\begin{align}
  u_g^{-1}=f_g \circ \exp\!.
\end{align}
Furthermore, let $s_g \colon U_g \rarr G$ be a local section of the principal $B$-bundle $p \colon G \rarr X$ defined by
\begin{align}
  s_g(x)=gf_g^{-1}(x)=gs_e(g^{-1}\cdot x)
\end{align}
for $x \in U_g$ and $g\in G$. The local section $s_g \colon U_g \rarr G$ gives a trivialization of the principal $B$-bundle $p\colon G \rarr X$ over $U_g$, i.e.\ $p^{-1}(U_g) \simeq U_g \times B$ as principal $B$-bundles over $U_g$.

For $\lambda \in \Hom_B(\mfrak{b},\C)$ there exists a $G$-equivariant sheaf of rings of twisted differential operators $\mcal{D}_X(\lambda)$ on $X$ such that $\mcal{D}_X(0) \simeq \mcal{D}_X$, where $\mcal{D}_X$ is the sheaf of rings of differential operators on $X$. If we denote by $\mcal{D}_X(\lambda)^{\rm op}$ the sheaf of rings opposite to $\mcal{D}_X(\lambda)$, then $\mcal{D}_X(\lambda)^{\rm op}$ is also a $G$-equivariant sheaf of rings of twisted differential operators on $X$, and we have
\begin{align}
  \mcal{D}_X(\lambda)^{{\rm op}} \simeq \mcal{D}_X(2\rho-\lambda), \label{eq:twisted sheaf op}
\end{align}
where the character $\rho \in \Hom_B(\mfrak{b},\C)$ of $\mfrak{b}$ is defined by
\begin{align}
  \rho(a) = -{\textstyle {1 \over 2}} \tr_{\mfrak{g}/\mfrak{b}} \ad(a) \label{eq:rho vector}
\end{align}
for $a \in \mfrak{b}$. Since $\mcal{D}_X(\lambda)$ is a $G$-equivariant sheaf of rings of twisted differential operators on $X$, we obtain a Lie algebra homomorphism
\begin{align}
  \alpha_{\mcal{D}_X(\lambda)}\colon \mfrak{g} \rarr \Gamma(X,\mcal{D}_X(\lambda)). \label{eq:Lie algebra embedding}
\end{align}
Now, let us consider a left $\mcal{D}_X(\lambda)$-module $\mcal{M}$. Then the vector space of global sections $\Gamma(X,\mcal{M})$ is a left $\Gamma(X,\mcal{D}_X(\lambda))$-module. Taking into account
\eqref{eq:Lie algebra embedding}, we get a $\mfrak{g}$-module structure on $\Gamma(X,\mcal{M})$.
A local section $s_g \colon U_g \rarr G$ gives an isomorphism
\begin{align}
  j^\lambda_{s_g} \colon \mcal{D}_X(\lambda)|_{U_g} \riso \mcal{D}_X|_{U_g}
\end{align}
of sheaves of rings of twisted differential operators on $U_g$. Hence, we define
\begin{align}
  \pi_\lambda^g(a)=(j^\lambda_{s_g} \circ \alpha_{\mcal{D}_X(\lambda)})(a)
\end{align}
for $a \in \mfrak{g}$.

Let $\{f_1,f_2,\dots,f_{\dim \widebar{\mfrak{n}}}\}$ be a basis of $\widebar{\mfrak{n}}$.
Then
\begin{align}
  u_g(x)=\sum_{i=1}^{\dim \widebar{\mfrak{n}}} u_g^i(x)f_i,\quad x \in U_g,  \label{eq:coordinates flag manifold}
\end{align}
where the functions
$u_g^i \colon U_g \rarr \C$ are called the coordinate functions on $U_g$.
\medskip

\theorem{\label{thm:operator realization}Let $\lambda \in \Hom_B(\mfrak{b},\C)$. Then we have
\begin{multline}
\pi^g_\lambda(a)= -\sum_{i=1}^{\dim \widebar{\mfrak{n}}} \bigg[{\ad(u_g(x))e^{\ad(u_g(x))} \over e^{\ad(u_g(x))}-{\rm id}_{\widebar{\mfrak{n}}}}\,(e^{-\ad(u_g(x))}\!\Ad(g^{-1})a)_{\widebar{\mfrak{n}}}\bigg]_i \partial_{u_g^i} \\
+\lambda((e^{-\ad(u_g(x))}\!\Ad(g^{-1})a)_\mfrak{b})
\end{multline}
for all $a\in \mfrak{g}$, where $x \in U_g$ and $[b]_i$ denotes the $i$-th coordinate of $b \in \widebar{\mfrak{n}}$ with respect to the basis $\{f_1,f_2,\dots,f_{\dim \widebar{\mfrak{n}}}\}$ of $\widebar{\mfrak{n}}$. In particular, we have
\begin{align}
    \pi^e_\lambda(a) = -\sum_{i=1}^{\dim \widebar{\mfrak{n}}} \bigg[{\ad(u_e(x)) \over e^{\ad(u_e(x))}-{\rm id}_{\widebar{\mfrak{n}}}}\,a\bigg]_i \partial_{u_e^i} \quad \text{for}\quad a \in \widebar{\mfrak{n}}
\end{align}
and
\begin{align}
    \pi^e_\lambda(a) = \sum_{i=1}^{\dim \widebar{\mfrak{n}}} [\ad(u_e(x))a]_i \partial_{u_e^i} + \lambda(a) \quad \text{for}\quad a \in \mfrak{h}.
\end{align}}

\proof{The proof of the parallel claim for any finite-dimensional semisimple Lie algebra and its parabolic subalgebra was given in \cite{Krizka-Somberg2015}. However, the assumption of the semisimplicity of $\mfrak{g}$ was redundant there, the only important structural data is a triangular decomposition of the Lie algebra in question. In particular, for the
 triangular
decomposition \eqref{eq:triangular decomposition ex} of our finite-dimensional (non-semisimple) Lie algebra $\mfrak{g}$ the proof given in \cite[Theorem 1.3]{Krizka-Somberg2015} is applicable verbatim to the present case.}

Let $\lambda \in \Hom_B(\mfrak{b},\C)$ be a character of $\mfrak{b}$. Then we denote by $\C_\lambda$ the $1$-dimensional representation of $\mfrak{b}$ defined by
\begin{align}
  av=\lambda(a)v \label{eq:1-dim representation}
\end{align}
for $a \in \mfrak{b}$ and $v \in \C_\lambda \simeq \C$ as vector spaces.
\medskip

\definition{Let $\lambda \in \Hom_B(\mfrak{b},\C)$ be a character of $\mfrak{b}$. Then the Verma module is the induced module
\begin{align}
  M^\mfrak{g}_\mfrak{b}(\lambda) = \Ind^\mfrak{g}_\mfrak{b} \C_\lambda \equiv U(\mfrak{g}) \otimes_{U(\mfrak{b})} \C_\lambda \simeq U(\widebar{\mfrak{n}}) \otimes_\C \C_\lambda,
\end{align}
where the last isomorphism of $U(\widebar{\mfrak{n}})$-modules follows from the Poincaré–Birkhoff–Witt
(PBW for short) theorem.}

\definition{Let $\lambda \in \Hom_B(\mfrak{b},\C)$ be a character of $\mfrak{b}$. Then a vector $v \in M^\mfrak{g}_\mfrak{b}(\lambda)$ is called a singular vector in the Verma module $M^\mfrak{g}_\mfrak{b}(\lambda)$ provided there exists a character $\mu \in \Hom_B(\mfrak{b},\C)$ of $\mfrak{b}$ such that $av=\mu(a)v$ for all $a \in \mfrak{b}$. We define a $\mfrak{b}$-module
\begin{align}
  M^\mfrak{g}_\mfrak{b}(\lambda)^\mfrak{b} = \langle \{v\in M^\mfrak{g}_\mfrak{b}(\lambda);\, av=\mu(a)v\ \text{for all}\ a \in \mfrak{b}\ \text{and some}\ \mu \in \Hom_B(\mfrak{b},\C)\} \rangle
\end{align}
and call it the vector space of singular vectors. We observe that singular vector $v \in M^\mfrak{g}_\mfrak{b}(\lambda)$ with character $\mu \in \Hom_B(\mfrak{b},\C)$ is also a weight vector with weight $\mu_{|\mfrak{h}} \in \mfrak{h}^*$.}

Let us consider a homomorphism
\begin{align}
  \varphi \colon M^\mfrak{g}_\mfrak{b}(\mu) \rarr M^\mfrak{g}_\mfrak{b}(\lambda)
\end{align}
of Verma modules. This homomorphism of $\mfrak{g}$-modules is uniquely determined by a homomorphism $\varphi_0 \colon \C_\mu \rarr M^\mfrak{g}_\mfrak{b}(\lambda)$ of $\mfrak{b}$-modules via the formula
\begin{align}
  \varphi(a \otimes v) = a\varphi_0(v)
\end{align}
for all $a \in U(\mfrak{g})$ and $v \in \C_\mu$. However, since $a\varphi_0(v)=\mu(a)\varphi_0(v)$ for all $a \in \mfrak{b}$ and $v \in \C_\mu$, we have $\varphi_0(\C_\mu) \subset M^\mfrak{g}_\mfrak{b}(\lambda)^\mfrak{b}$. Therefore, we obtain a homomorphism
\begin{align}
  \varphi_0 \colon \C_\mu \rarr M^\mfrak{g}_\mfrak{b}(\lambda)^\mfrak{b}
\end{align}
of $\mfrak{b}$-modules, which uniquely determines $\varphi$. In the opposite way, any such homomorphism $\varphi_0$ gives rise to a homomorphism $\varphi \colon M^\mfrak{g}_\mfrak{b}(\mu) \rarr M^\mfrak{g}_\mfrak{b}(\lambda)$ of Verma modules. Therefore, this gives us a vector space isomorphism
\begin{align}
  \Hom_\mfrak{g}(M^\mfrak{g}_\mfrak{b}(\mu), M^\mfrak{g}_\mfrak{b}(\lambda))  \simeq \Hom_\mfrak{b}(\C_\mu, M^\mfrak{g}_\mfrak{b}(\lambda)^\mfrak{b}) .
\end{align}

Let $\{f_1,f_2,\dots,f_{\dim \widebar{\mfrak{n}}}\}$ be a basis of $\widebar{\mfrak{n}}$. We denote by $\{x_1,x_2,\dots,x_{\dim \widebar{\mfrak{n}}}\}$ the linear coordinate functions on $\widebar{\mfrak{n}}$ with respect to the given basis of $\widebar{\mfrak{n}}$ and by $\{y_1,y_2,\dots,y_{\dim \widebar{\mfrak{n}}}\}$ the dual linear coordinate functions on $\widebar{\mfrak{n}}$. The Weyl algebra $\eus{A}^\mfrak{g}_{\widebar{\mfrak{n}}}$ of the vector space $\widebar{\mfrak{n}}$ is generated by $\{x_i, \partial_{x_i};\, i=1,2,\dots, \dim \widebar{\mfrak{n}}\}$ and the Weyl algebra $\eus{A}^\mfrak{g}_{\widebar{\mfrak{n}}^*}$ of the vector space $\widebar{\mfrak{n}}^*$ by $\{y_i, \partial_{y_i};\, i=1,2,\dots, \dim \widebar{\mfrak{n}}\}$. Furthermore, we define the algebraic Fourier transform
\begin{align}
  \mcal{F} \colon \eus{A}^\mfrak{g}_{\widebar{\mfrak{n}}} \riso  \eus{A}^\mfrak{g}_{\widebar{\mfrak{n}}^*}
\end{align}
as an isomorphism of $\C$-algebras uniquely determined by
\begin{align}
  \mcal{F}(x_i)= -\partial_{y_i}, \qquad \mcal{F}(\partial_{x_i}) = y_i
\end{align}
for all $i=1,2,\dots, \dim \widebar{\mfrak{n}}$. We denote by $\eus{I}_e$ the left ideal of $\eus{A}^\mfrak{g}_{\widebar{\mfrak{n}}}$ generated by polynomials on $\widebar{\mfrak{n}}$ vanishing at the point $0$.

Since the mapping $u_e \colon U_e \rarr \widebar{\mfrak{n}}$ is a diffeomorphism, it induces an isomorphism
\begin{align}
  \Psi_{u_e} \colon  \mcal{D}_X({U_e}) \riso \Gamma(\widebar{\mfrak{n}},\mcal{D}_{\widebar{\mfrak{n}}})
\end{align}
of associative $\C$-algebras. Let us note that the Weyl algebra $\eus{A}^\mfrak{g}_{\widebar{\mfrak{n}}}$ is contained in $\Gamma(\widebar{\mfrak{n}},\mcal{D}_{\widebar{\mfrak{n}}})$. Moreover, we have
\begin{align}
  \Psi_{u_e}(u_e^i)=x_i \quad \text{and} \quad \Psi_{u_e}(\partial_{u_e^i})=\partial_{x_i}  \label{eq:Psi transformation}
\end{align}
for $i=1,2,\dots, \dim \widebar{\mfrak{n}}$, where the coordinate functions $\{u_e^i;\, i=1,2,\dots, \dim \widebar{\mfrak{n}}\}$ on $U_e$ are given by \eqref{eq:coordinates flag manifold}. By a completely analogous way as in \cite[Section 2]{Krizka-Somberg2015}, we obtain an isomorphism of $\mfrak{g}$-modules
\begin{align}
  \Phi_\lambda \colon M^\mfrak{g}_\mfrak{b}(\lambda-\rho) \riso \eus{A}^\mfrak{g}_{\widebar{\mfrak{n}}}/\eus{I}_e
\end{align}
uniquely determined by
\begin{align}
  \Phi_\lambda(v_{\lambda-\rho})= 1\ {\rm mod}\ \eus{I}_e,
\end{align}
where $v_{\lambda-\rho}$ is the highest weight vector of $M^\mfrak{g}_\mfrak{b}(\lambda-\rho)$. The structure of a left $\mfrak{g}$-module on $\eus{A}^\mfrak{g}_{\widebar{\mfrak{n}}}/\eus{I}_e$ is given through a homomorphism
\begin{align}
  \pi_\lambda \colon \mfrak{g} \rarr \eus{A}^\mfrak{g}_{\widebar{\mfrak{n}}}
\end{align}
of Lie algebras defined by
\begin{align}
  \pi_\lambda = \Psi_{u_e} \circ \pi^e_{\lambda+\rho}. \label{eq:action pi_lambda}
\end{align}
In the next theorem we realize the mapping $\Phi_\lambda$ in an explicit way. This result enables us to write singular vectors as elements of $M^\mfrak{g}_\mfrak{b}(\lambda-\rho) \simeq U(\widebar{\mfrak{n}}) \otimes_\C \C_{\lambda-\rho}$.
\medskip

\theorem{\label{thm:symmetrization}
Let $\{f_1,f_2,\dots,f_{\dim \widebar{\mfrak{n}}}\}$ be a basis of $\widebar{\mfrak{n}}$ and $\{x_1,x_2,\dots,x_{\dim \widebar{\mfrak{n}}}\}$ be the corresponding linear coordinate functions on $\widebar{\mfrak{n}}$. Further, let $\beta \colon S(\widebar{\mfrak{n}}) \rarr U(\widebar{\mfrak{n}})$ be the symmetrization map defined by
\begin{align}
  \beta(a_1a_2\dots a_k) = {1 \over k!} \sum_{\sigma \in S_k} a_{\sigma(1)} a_{\sigma(2)} \dots a_{\sigma(k)}
\end{align}
for all $k \in \N$ and $a_1,a_2,\dots, a_k \in \widebar{\mfrak{n}}$. Then
\begin{align}
  \Phi_\lambda(\beta(f_{i_1}f_{i_2} \dots f_{i_k})v_{\lambda-\rho})= (-1)^k \partial_{x_{i_1}} \partial_{x_{i_2}} \dots \partial_{x_{i_k}}\ {\rm mod}\ \eus{I}_e
\end{align}
for all $k \in \N$ and $i_1,i_2,\dots,i_k \in \{1,2,\dots,\dim \widebar{\mfrak{n}}\}$.}

\proof{The proof goes along the same line as that in \cite[Theorem 2.4]{Krizka-Somberg2015}.}

The algebraic Fourier transform $\mcal{F} \colon \eus{A}^\mfrak{g}_{\widebar{\mfrak{n}}} \rarr \eus{A}^\mfrak{g}_{\widebar{\mfrak{n}}^*}$ leads to an isomorphism
\begin{align}
  \tau \colon \eus{A}^\mfrak{g}_{\widebar{\mfrak{n}}}/\eus{I}_e \riso \eus{A}^\mfrak{g}_{\widebar{\mfrak{n}}^*}/\mcal{F}(\eus{I}_e)
\end{align}
of $\mfrak{g}$-modules defined by
\begin{align}
  Q\ {\rm mod}\ \eus{I}_e \mapsto \mcal{F}(Q)\ {\rm mod}\ \mcal{F}(\eus{I}_e)
\end{align}
for all $Q \in \eus{A}^\mfrak{g}_{\widebar{\mfrak{n}}}$. The structure of a left $\mfrak{g}$-module on $\eus{A}^\mfrak{g}_{\widebar{\mfrak{n}}^*}/\mcal{F}(\eus{I}_e)$ is given through a homomorphism \begin{align}
  \hat{\pi}_\lambda \colon \mfrak{g} \rarr \eus{A}^\mfrak{g}_{\widebar{\mfrak{n}}^*}
\end{align}
of Lie algebras, defined by
\begin{align}
  \hat{\pi}_\lambda = \mcal{F} \circ \pi_\lambda.
\end{align}
Since there is a canonical isomorphism of left $\eus{A}^\mfrak{g}_{\widebar{\mfrak{n}}^*}$-modules
\begin{align}
  \C[\widebar{\mfrak{n}}^*] \riso \eus{A}^\mfrak{g}_{\widebar{\mfrak{n}}^*}/ \mcal{F}(\eus{I}_e),
\end{align}
we obtain immediately an isomorphism
\begin{align}
  \tau \circ \Phi_\lambda \colon M^\mfrak{g}_\mfrak{b}(\lambda-\rho) \riso \C[\widebar{\mfrak{n}}^*] \label{eq:Verma module isomorphism}
\end{align}
of $\mfrak{g}$-modules. If we introduce a $\mfrak{b}$-module
\begin{align}
  \Sol(\mfrak{g},\mfrak{b}; \C[\widebar{\mfrak{n}}^*])^\mcal{F} =\{u \in \C[\widebar{\mfrak{n}}^*];\, \hat{\pi}_\lambda(a)u=\mu(a)u\ \text{for all}\ a \in \mfrak{b}\ \text{and some}\ \mu \in \Hom_B(\mfrak{b},\C)\},
\end{align}
then by \eqref{eq:Verma module isomorphism}, we obtain an isomorphism of $\mfrak{b}$-modules
\begin{align}
  \tau \circ \Phi_\lambda \colon M^\mfrak{g}_\mfrak{b}(\lambda-\rho)^\mfrak{b} \riso \Sol(\mfrak{g},\mfrak{b}; \C[\widebar{\mfrak{n}}^*])^\mcal{F} .
\end{align}

Therefore, the algebraic problem of finding singular vectors in the Verma module $M^\mfrak{g}_\mfrak{b}(\lambda-\rho)$ is converted into an analytic problem of solving a system of partial differential equations on $\C[\widebar{\mfrak{n}}^*]$.


\section{Conformal Galilei algebras $\mathfrak{cga}_\ell(d,\R)$ and their representation theory}
\label{sec:confgalalg}

Conformal Galilei algebras $\mfrak{cga}_\ell(d,\R)$ are generalized non-relativistic versions of conformal Lie algebras $\mfrak{so}(d+1,2,\R)$ and parameterized by pairs $(d,\ell)$ with
$d\in\mN$ (called the spatial dimension) and $\ell \in {1 \over 2}\N_0$ (called the spin). The Lie algebras $\mfrak{cga}_\ell(d,\R)$ are not semisimple, and there exists a $1$-dimensional central extension either for $d\in \N$ and $\ell$ half-integer or for $d=2$ and $\ell$ non-negative integer.

Let $\R^{1,d}$ be the $(d+1)$-dimensional space endowed with the canonical metric of signature $(1,d)$ and let $(t,x_1,x_2,\dots,x_d)$ be the canonical linear coordinate functions on $\R^{1,d}$.  The variable $t$ is referred to as the time coordinate while $x_1,x_2,\dots,x_d$ are referred to as the spatial coordinates. Let us consider the infinitesimal transformations of $\R^{1,d}$ given by vector fields
\begin{align}\label{eq:cga relations}
\begin{gathered}
  H=\partial_t, \qquad D=-2t\partial_t-2\ell\sum_{i=1}^d x_i \partial_{x_i}, \qquad C=t^2\partial_t + 2\ell \sum_{i=1}^d tx_i\partial_{x_i} \\
  M_{ij}=-x_i\partial_{x_j}+x_j \partial_{x_i}, \qquad P_{n,i}=(-t)^n\partial_{x_i},
\end{gathered}
\end{align}
where $i,j=1,2\dots,d$ and $n=0,1,\dots,2\ell$ with $\ell \in {1 \over 2}\N_0$. Their linear span has real
dimension $\smash{{d(d-1) \over 2}}+(2\ell+1)d + 3$ and is closed under the Lie bracket of vector fields. Consequently, we get finite-dimensional real Lie algebra with non-trivial Lie brackets
\begin{align}
\begin{gathered}[]
  [D,H]=2H, \qquad \qquad [C,H]=D, \qquad \qquad [D,C]=-2C,\\
  [H,P_{n,i}]=-nP_{n-1,i}, \qquad [D,P_{n,i}]=2(\ell-n)P_{n,i}, \qquad [C,P_{n,i}]=(2\ell-n)P_{n+1,i},\\
  [M_{ij},M_{k\ell}]=-\delta_{i,k}M_{j\ell}-\delta_{j,\ell}M_{ik}+ \delta_{i,\ell}M_{jk} + \delta_{j,k}M_{i\ell}, \\
  [M_{ij},P_{n,k}]=-\delta_{i,k}P_{n,j}+\delta_{j,k}P_{n,i},
\end{gathered}
\end{align}
where $i,j=1,2\dots,d$ and $n=0,1,\dots,2\ell$. The Lie algebra $\mfrak{cga}_\ell(d,\R)$ has a Lie subalgebra given by direct sum of $\mfrak{sl}(2,\R) \simeq \mfrak{so}(2,1,\R)$ generated by $C,D,H$ and $\mfrak{so}(d,\R)$ generated by $M_{ij}$ for $i,j=1,2,\dots,d$. The subspace generated by $P_{n,i}$ for $i=1,2,\dots,d$ and $n=0,1,\dots,2\ell$ forms an abelian ideal and we have $\mfrak{cga}_\ell(d,\R) \simeq (\mfrak{sl}(2,\R) \oplus \mfrak{so}(d,\R)) \ltimes \mathbb{V}_{2\ell \omega,\omega_1}$, where $\mathbb{V}_{2\ell \omega,\omega_1}$ is the irreducible finite-dimensional $(\mfrak{sl}(2,\R) \oplus \mfrak{so}(d,\R))$-module with highest weight $(2\ell \omega,\omega_1)$ (see the next section for
the representation theoretical conventions concerning weights $\omega,\omega_1$).

The conformal Galilei algebras $\mfrak{cga}_\ell(d,\R)$ admit two distinct types of central extensions according to the values of $d\in \N$ and $\ell\in {1 \over 2}\N_0$. In particular, we have
\begin{enum}
  \item[1)] the mass central extension for $d \in \N$ and $\ell \in {1\over 2}+\N_0$,
  \begin{align}
    [P_{m,i},P_{n,j}]=\delta_{i,j}\delta_{m+n,2\ell}I_m M, \qquad I_m=(-1)^{m+\ell+{1\over 2}}(2\ell-m)!m!,
  \end{align}
  \item[2)] the exotic central extension for $d=2$ and $\ell \in \N_0$,
  \begin{align}
    [P_{m,i},P_{n,j}]=\veps_{i,j}\delta_{m+n,2\ell}I_m \Theta, \qquad I_m=(-1)^m(2\ell-m)!m!,
  \end{align}
  where $\veps_{1,2}=-\veps_{2,1}=1$ and $\veps_{1,1}=\veps_{2,2}=0$.
\end{enum}

Let us note that the generators of time translation $H$, space translations $P_{0,i}$, spatial rotations $M_{ij}$ and Galilei transformations $P_{1,i}$ form the Galilei algebra of $\R^{1,d}$.
\medskip

In the rest of the article we shall restrict to the (complexification of the real) conformal Galilei algebras
$\mfrak{cga}_\ell(d,\C)$ for $d=1$ and $\ell \in \N$, and we denote $P_{n,0}$ by $P_n$ for $n=0,1,\dots,2\ell$. In particular, as follows from the previous paragraph there is no central extension in this case.


\subsection{Representation theoretical conventions}

The Lie algebra $\mathfrak{cga}_\ell(1,\C)$ is given by the complex vector space
$\langle D, H, C, P_0, P_1, \dots , P_{2\ell}\rangle$ together with the following non-trivial Lie brackets
\begin{align}\label{confgalalg}
\begin{aligned}
& [D,H]=2H, &\quad & [C,H]=D, &\quad &  [D,C]=-2C, \\
& [H,P_n]=-nP_{n-1}, & & [D,P_n]=2(\ell-n)P_n, & &  [C,P_n]=(2\ell-n)P_{n+1}
\end{aligned}
\end{align}
for all $n=0,1,\dots ,2\ell$.
\medskip

\remark If we denote
\begin{align}
\begin{gathered}
  L_{-1}=H,\qquad L_0={\textstyle {1 \over 2}}D,\qquad L_1=C, \\
  M_n=P_{n+\ell}
\end{gathered}
\end{align}
for $n=-\ell,-\ell+1,\dots,\ell$, then we may write
\begin{align}
  [L_m,L_n]=(m-n)L_{m+n}, \qquad [L_m,M_n]=(\ell m-n)M_{m+n}, \qquad [M_m,M_n]=0. \label{eq:infinite extension}
\end{align}
Now, for $\ell \in \N_0$ we may define an infinite-dimensional complex Lie algebra with a basis $L_m,M_n$ for $m,n\in \Z$ together with the Lie bracket given by \eqref{eq:infinite extension} for all $m,n\in \Z$. This infinite-dimensional complex Lie algebra admits a universal central extension whose explicit form depends on the value of $\ell \in \N_0$. For $\ell=0$ we obtain the so called Heisenberg-Virasoro algebra and for $\ell=1$ we get the so called Bondi-Metzner-Sachs (BMS) algebra in dimension $3$. This observation was the first evidence for the non-relativistic version of ${\rm AdS}_3/{\rm CFT}_2$ correspondence, see \cite{Bagchi-Gopakumar2009}, \cite{Martelli-Tachikawa2010}.
\medskip

Let us denote $\mfrak{g}=\mfrak{cga}_\ell(1,\C)$, $\ell \in \N$. Let $G$ be a connected complex Lie group with its Lie algebra $\mfrak{g}$. By structural theory reviewed in Section \ref{sec:geomrealindmod}, we choose the Cartan subalgebra $\mfrak{h}$ of $\mfrak{g}$ to be
\begin{align}
  \mfrak{h} = \C D \oplus \C P_\ell.
\end{align}
We define $\omega_D,\omega_{P_\ell} \in \mfrak{h}^*$ by
\begin{align}
\omega_D(D)=1,\quad \omega_D(P_\ell)=0,\quad \omega_{P_\ell}(D)=0,\quad \omega_{P_\ell}(P_\ell)=1.
\end{align}
Then the root system of $\mfrak{g}$ with respect to $\mfrak{h}$ is $\Delta=\{\pm k\alpha;\,1\leq k \leq \ell\}$, $\alpha=2\omega_D$, and the root spaces are
\begin{align}
  \mfrak{g}_{k\alpha}= \C P_{\ell-k}, \qquad \mfrak{g}_\alpha = \C H \oplus \C P_{\ell-1},\qquad  \mfrak{g}_{-\alpha} = \C C \oplus \C P_{\ell+1}, \qquad \mfrak{g}_{-k\alpha}= \C P_{\ell+k}
\end{align}
for $k=2,3,\dots,\ell$. Further, the positive root system is $\Delta^+=\{k\alpha;\, 1\leq k \leq \ell\}$ and the negative root system is $\Delta^-=\{-k\alpha;\, 1\leq k \leq \ell\}$.

The Lie algebra $\mfrak{g}$ is $|\ell|$-graded with respect to the grading given by $\mfrak{g}_i= \mfrak{g}_{i\alpha}$ for $0\neq i \in \Z$ and $\mfrak{g}_0=\mfrak{h}$. We define nilpotent Lie subalgebras $\mfrak{n}_+$ and $\mfrak{n}_-$ of $\mfrak{g}$ by
\begin{align}
 \mathfrak{n}_+=\bigoplus_{i=1}^\ell \mfrak{g}_i= \langle H, P_0, P_1, \dots , P_{\ell-1}\rangle,\qquad \mathfrak{n}_-=\bigoplus_{i=1}^\ell \mfrak{g}_{-i} =\langle C, P_{\ell+1}, P_{\ell+2},\dots ,P_{2\ell}\rangle,
\end{align}
and the standard Borel subalgebra $\mfrak{b}$ of $\mfrak{g}$ by
\begin{align}
\mfrak{b}= \mfrak{h} \oplus \mfrak{n}_-=\langle D,C,P_\ell,P_{\ell+1},\dots,P_{2\ell} \rangle.
\end{align}
We remark that $\mfrak{b}$ is contained in a maximal (standard) solvable Lie subalgebra of $\mfrak{g}$, which is given by the linear span $\langle D,C,P_0,P_1,\dots,P_{2\ell} \rangle$. Moreover, we have a triangular decomposition
\begin{align}\label{triangdecomp}
\mfrak{g} = \mfrak{n}_-\oplus \mfrak{h} \oplus \mfrak{n}_+
\end{align}
of the Lie algebra $\mfrak{g}$ with $[\mathfrak{h}, \mathfrak{n}_+]= \mathfrak{n}_+$ and $[\mathfrak{h}, \mathfrak{n}_-]= \mathfrak{n}_-$. The Borel subgroup $B$ of $G$ with its Lie algebra $\mfrak{b}$ is defined by $B=N_G(\mfrak{b})$.

Any character $\lambda \in \Hom_B(\mfrak{b},\C)$ of $\mfrak{b}$ is given by
\begin{align}
  \lambda=\delta \widetilde{\omega}_D + p \widetilde{\omega}_{P_\ell}
\end{align}
for some $\delta, p \in \C$, where $\widetilde{\omega}_D, \widetilde{\omega}_{P_\ell} \in \Hom_B(\mfrak{b},\C)$ are equal to $\omega_D, \omega_{P_\ell} \in \mfrak{h}^*$ regarded as trivially extended to $\mfrak{b}=\mfrak{h} \oplus \mfrak{n}_-$. We denote the character $\delta \widetilde{\omega}_D+p \widetilde{\omega}_{P_\ell}$ for $\delta,p \in \C$ by $\lambda_{\delta,p}$. The character $\rho \in \Hom_B(\mfrak{b},\C)$ of $\mfrak{b}$ introduced in \eqref{eq:rho vector} is
\begin{align}
\rho=-\big(1+{\textstyle {1\over 2}}\ell(\ell+1)\!\big)\widetilde{\omega}_D.
\end{align}
We shall denote
\begin{align}
\delta_\rho =-\big(1+{\textstyle {1\over 2}}\ell(\ell+1)\!\big) \qquad \text{and} \qquad p_\rho=0.
\end{align}

Since for $\delta,p \in \C$ the character $\lambda_{\delta,p} \in \Hom_B(\mfrak{b},\C)$ of $\mfrak{b}$ defines the $1$-dimensional representation $\C_{\delta,p}$ of $\mfrak{b}$ through the formula \eqref{eq:1-dim representation}, we associate to $\C_{\delta,p}$ the Verma module $M_\ell(\delta,p)$:
\begin{align}\label{eq:verma module}
M_\ell(\delta,p)=U(\mathfrak{g})\otimes_{U(\mathfrak{b})}\mC_{\delta,p}
\simeq U(\mathfrak{n}_+)\otimes_\C\mC_{\delta,p},
\end{align}
where the isomorphism of $U(\mathfrak{n}_+)$-modules follows from the PBW theorem. We denote by $v_{\delta,p}$ the lowest weight vector of $M_\ell(\delta,p)$ with lowest weight $\delta \omega_D+ p\omega_{P_\ell}$.
There is an automorphism $\omega \colon \mathfrak{g} \rarr \mathfrak{g}$ of $\mfrak{g}$ defined by
\begin{align}
\omega(D)=-D,\quad \omega(H)=-C,\quad \omega(C)=-H,\quad \omega(P_n)=-P_{2\ell-n}
\end{align}
for $n=0,1,\dots,2\ell$, called the Chevalley involution. Moreover, there exists a unique symmetric bilinear form $\langle \cdot\,,\cdot \rangle$ on $M_\ell(\delta,p)$, called the Shapovalov form, with its contravariance property
$\langle v_{\delta,p}, v_{\delta,p} \rangle =1$ and
$\langle au, v \rangle = -\langle u, \omega(a)v \rangle$ for all
$u,v\in M_\ell(\delta,p)$ and $a\in \mathfrak{g}$.


\subsection{Geometric realization of Verma modules}

Let us denote by $\{t,x_0,\dots,x_{\ell-1}\}$ the linear coordinate functions on $\mfrak{n}_+$ with respect to the given basis $\{H,P_0,\dots,P_{\ell-1}\}$ of $\mfrak{n}_+$, and by $\{u,y_0,\dots,y_{\ell-1}\}$ the dual linear coordinate functions on $\mfrak{n}_+^*$. Then the Weyl algebra $\eus{A}^\mfrak{g}_{\mfrak{n}_+}$ is generated by
\begin{align}
  \{t,x_0,\dots,x_{\ell-1},\partial_t,\partial_{x_0}, \dots, \partial_{x_{\ell-1}}\}
\end{align}
and the Weyl algebra $\smash{\eus{A}^\mfrak{g}_{\mfrak{n}_+^*}}$ by
\begin{align}
  \{u,y_0,\dots,y_{\ell-1},\partial_u,\partial_{y_0}, \dots, \partial_{y_{\ell-1}}\}.
\end{align}
Furthermore, the coordinate functions $\{u_e^t,u_e^{x_0},\dots,u_e^{x_{\ell-1}}\}$ on $U_e$ are defined by
\begin{align}
  u_e(x)=u^t_e(x)H+\sum_{n=0}^{\ell-1} u_e^{x_n}(x) P_n \label{eq:u_e(x) coordinates}
\end{align}
for $x \in U_e$. Using \eqref{eq:Psi transformation}, we obtain
\begin{align}
  \Psi_{u_e}\!(u_e^t)=t,\qquad \Psi_{u_e}\!(u_e^{x_n})=x_n, \qquad  \Psi_{u_e}\!(\partial_{u_e^t})=\partial_t,\qquad \Psi_{u_e}\!(\partial_{u_e^{x_n}})=\partial_{x_n}
\end{align}
for $n=0,1,\dots,\ell-1$. We define
\begin{align}
 u(x)= t(x) H+\sum_{n=0}^{\ell-1} x_n(x) P_n \label{eq:u(x) coordinates}
\end{align}
for $x \in \mfrak{n}_+$. Finally, the homomorphism $\pi_{\delta,p} \colon \mfrak{g} \rarr \eus{A}^\mfrak{g}_{\mfrak{n}_+}$ of Lie algebras defined by \eqref{eq:action pi_lambda} is
\begin{align}
  \pi_{\delta,p} = \Psi_{u_e} \circ \pi^e_{\lambda_{\delta,p} + \rho},
\end{align}
and its explicit form is a subject of the next theorem.
\medskip

Let us recall that the Bernoulli polynomials $B_n(x)$ are defined by a generating function
\begin{align}
  {t e^{xt} \over e^t-1} = \sum_{n=0}^\infty B_n(x) {t^n \over n!}  \label{eq:Bernoulli polynomials}
\end{align}
for $0 \neq t \in \R$ and $x \in \R$. We denote $B_n=B_n(0)$ the first Bernoulli numbers.
\medskip

\theorem{\label{thm:cga realization}Let $\delta, p \in \C$. Then the embedding of $\mfrak{g}$ into $\eus{A}^\mfrak{g}_{\mfrak{n}_+}$ and $\eus{A}^\mfrak{g}_{\mfrak{n}_+^*}$
is given by
\begin{enumerate}
  \item[1)]
  \begin{align}
  \begin{aligned}
    \pi_{\delta,p}(H)&= -\partial_t-  \sum_{j=1}^{\ell-1}\sum_{k=j}^{\ell-1} (-1)^{j-1}B_j \binom{k}{j}t^{j-1}x_k \partial_{x_{k-j}}, \\
    \pi_{\delta,p}(P_n)&=-\partial_{x_n}-\sum_{j=1}^n (-1)^jB_j\binom{n}{j}t^j \partial_{x_{n-j}},
  \end{aligned}
  \end{align}
  \begin{align}
  \begin{aligned}
    \hat{\pi}_{\delta,p}(H)&=-u+\sum_{j=1}^{\ell-1}\sum_{k=j}^{\ell-1} B_j\binom{k}{j}y_{k-j}\partial_{y_k}\partial_u^{j-1},\\
    \hat{\pi}_{\delta,p}(P_n)&= -y_n-\sum_{j=1}^n B_j\binom{n}{j}y_{n-j}\partial^j_u
  \end{aligned}
  \end{align}
  for $n=0,1,\dots,\ell-1$;
  \item[2)]
  \begin{align}
  \begin{aligned}
    \pi_{\delta,p}(P_\ell)&=-\ell t\partial_{x_{\ell-1}}+p, \\
    \pi_{\delta,p}(D)&=-2t\partial_t - 2\sum_{j=0}^{\ell-1} (\ell-j)x_j\partial_{x_j} +
    \delta+\delta_\rho,
  \end{aligned}
  \end{align}
  \begin{align}
  \begin{aligned}
    \hat{\pi}_{\delta,p}(P_\ell)&=\ell y_{\ell-1}\partial_u +p, \\
    \hat{\pi}_{\delta,p}(D)&=2u\partial_u + 2\sum_{j=0}^{\ell-1} (\ell-j)y_j\partial_{y_j} +
    \delta-\delta_\rho;
  \end{aligned}
  \end{align}
  \item[3)]
  \begin{align}
  \begin{aligned}
    \pi_{\delta,p}(C)&= -t^2\partial_t -\sum_{j=0}^{\ell-1} (\ell-j)tx_j \partial_{x_j} -\sum_{j=0}^{\ell-2} (2\ell-j)x_j\partial_{x_{j+1}} + (\delta+\delta_\rho)t \\ &\quad
    + \sum_{j=1}^\ell B_j(\ell+1) \binom{\ell}{j} t^j x_{\ell-1} \partial_{x_{\ell-j}} +p(\ell+1)x_{\ell-1}
    \\ & \quad + \sum_{j=2}^{\ell-1} \sum_{k=0}^{\ell-2} {2(\ell-1-k) \over k+1} B_j\binom{k+1}{j} t^j x_k \partial_{x_{k-j+1}}, \\
    \pi_{\delta,p}(P_n)&= -\sum_{j=0}^{\ell-1} \sum_{k=n-\ell+1}^{n-j} B_j \binom{n}{k}\!\binom{n-k}{j} t^{k+j} \partial_{x_{n-k-j}} +  p\binom{n}{\ell}t^{n-\ell}
  \end{aligned}
  \end{align}
  \begin{align}
  \begin{aligned}
    \hat{\pi}_{\delta,p}(C)&= -u\partial^2_u - \sum_{j=0}^{\ell-1} (\ell-j)y_j \partial_{y_j}\partial_u + \sum_{j=0}^{\ell-2} (2\ell-j)y_{j+1}\partial_{y_j} - (\delta-\delta_\rho)\partial_u  \\ &\quad
     -\sum_{j=1}^\ell (-1)^jB_j(\ell+1) \binom{\ell}{j} y_{\ell-j} \partial_{y_{\ell-1}}\partial^j_u -p(\ell+1)\partial_{y_{\ell-1}} \\ &\quad - \sum_{j=2}^{\ell-1} \sum_{k=0}^{\ell-2} (-1)^j{2(\ell-1-k) \over k+1} B_j\binom{k+1}{j} y_{k-j+1} \partial_{y_k}\partial_u^j, \\
    \hat{\pi}_{\delta,p}(P_n)&= -\sum_{j=0}^{\ell-1} \sum_{k=n-\ell+1}^{n-j} (-1)^{k+j}B_j \binom{n}{k}\!\binom{n-k}{j} y_{n-k-j} \partial^{k+j}_u +  (-1)^{n-\ell}p\binom{n}{\ell}\partial^{n-\ell}_u
  \end{aligned}
  \end{align}
  for $n=\ell+1,\ell+2,\dots,2\ell$.
\end{enumerate}}

\proof{In what follows we use the notation $\ad_{u(x)}$ instead of $\ad(u(x))$. Using \eqref{eq:u(x) coordinates}, we obtain
\begin{align} \label{eq:ad(u_e(x)) relation k=1}
\begin{aligned}
  \ad_{u(x)}(P_n)&={\textstyle \sum\limits_{j=0}^{\ell-1}}x_j[P_j,P_n]+t[H,P_n]=-ntP_{n-1}, \\
  \ad_{u(x)}(H)&={\textstyle \sum\limits_{j=0}^{\ell-1}}x_j[P_j,H]+t[H,H]={\textstyle \sum\limits_{j=0}^{\ell-1}} jx_jP_{j-1}, \\
  \ad_{u(x)}(D)&={\textstyle \sum\limits_{j=0}^{\ell-1}}x_j[P_j,D]+t[H,D]=-{\textstyle \sum\limits_{j=0}^{\ell-1}} 2(\ell-j)x_jP_j-2tH, \\
  \ad_{u(x)}(C)&={\textstyle \sum\limits_{j=0}^{\ell-1}}x_j[P_j,C]+t[H,C]=-{\textstyle \sum\limits_{j=0}^{\ell-1}} (2\ell-j)x_jP_{j+1}-tD.
\end{aligned}
\end{align}
By a similar computation to above, we find
\begin{align} \label{eq:ad(u_e(x)) relation k>1}
\begin{aligned}
  \ad^k_{u(x)}(P_n)&=(-1)^k k! \binom{n}{k} t^kP_{n-k}, \\
  \ad^k_{u(x)}(H)&={\textstyle \sum\limits_{j=0}^{\ell-1}} (-1)^{k-1} k! \binom{j}{k} t^{k-1}x_jP_{j-k}
\end{aligned}
\end{align}
for all $k \in \N$, and
\begin{align} \label{eq:ad(u_e(x)) relation C,D}
\begin{aligned}
  \ad^2_{u(x)}(C)&= {\textstyle \sum\limits_{j=0}^{\ell-1}} \big((2\ell-j)(j+1)+2(\ell-j)\big)tx_jP_j + 2t^2H, \\
  \ad^k_{u(x)}(D)&= {\textstyle \sum\limits_{j=0}^{\ell-2}} 2(-1)^k{\ell-1-j \over j+1} k!\binom{j+1}{k} t^{k-1}x_jP_{j-k+1}
\end{aligned}
\end{align}
for $k \geq 2$. By Theorem \ref{thm:operator realization}, we have
\begin{align*}
  \pi_{\delta,p}(a)= \sum_{i=1}^{\ell+1}\,[\ad_{u(x)}(a)]_i\partial_{u^i}+(\lambda_{\delta,p}+\rho)(a)
\end{align*}
for $a \in \mfrak{h}$, where $u^i=x_{i-1}$ for $i=1,2,\dots,\ell$ and $u^{\ell+1}=t$. Therefore, using \eqref{eq:ad(u_e(x)) relation k=1}, we get
\begin{align*}
  \pi_{\delta,p}(D)&=-2t\partial_t-2\sum_{j=0}^{\ell-1}(\ell-j)x_j\partial_{x_j}+\delta-\big(1+{\textstyle {1\over 2}}\ell(\ell+1)\!\big), \\
  \pi_{\delta,p}(P_\ell)&=-\ell t\partial_{x_{\ell-1}}+p.
\end{align*}
Similarly, from Theorem \ref{thm:operator realization} we have
\begin{align*}
  \pi_{\delta,p}(a)= -\sum_{i=1}^{\ell+1}\bigg[{\ad_{u(x)} \over e^{\ad_{u(x)}}-{\rm id}_{\mfrak{n}_+}}\,a\bigg]_i\partial_{u^i}
\end{align*}
for $a \in \mfrak{n}_+$, where $u^i=x_{i-1}$ for $i=1,2,\dots,\ell$ and $u^{\ell+1}=t$, which we can rewrite using \eqref{eq:Bernoulli polynomials} into the form
\begin{align*}
  \pi_{\delta,p}(a)=-\sum_{i=1}^{\ell+1} \sum_{j=0}^\infty {B_j \over j!}\, [\ad^j_{u(x)}(a)]_i \partial_{u^i}
\end{align*}
for $a \in \mfrak{n}_+$. Therefore, we get
\begin{align*}
  \pi_{\delta,p}(P_n)&= -\sum_{i=1}^{\ell+1} \sum_{j=0}^\infty {B_j \over j!}\, [\ad^j_{u(x)}(P_n)]_i \partial_{u^i}
  = -\partial_{x_n} - \sum_{j=1}^\infty (-1)^j B_j \binom{n}{j} t^j \partial_{x_{n-j}} \\
  & = -\partial_{x_n} - \sum_{j=1}^n (-1)^j B_j \binom{n}{j} t^j \partial_{x_{n-j}}, \\
  \pi_{\delta,p}(H)&=-\sum_{i=1}^{\ell+1} \sum_{j=0}^\infty {B_j \over j!}\, [\ad^j_{u(x)}(H)]_i \partial_{u^i} = -\partial_t - \sum_{j=1}^\infty \sum_{k=0}^{\ell-1} (-1)^{j-1}B_j \binom{k}{j}t^{j-1}x_k \partial_{x_{k-j}} \\
  &= -\partial_t - \sum_{j=1}^{\ell-1} \sum_{k=j}^{\ell-1} (-1)^{j-1}B_j \binom{k}{j}t^{j-1}x_k \partial_{x_{k-j}},
\end{align*}
where we used \eqref{eq:ad(u_e(x)) relation k>1}.

Finally, from Theorem \ref{thm:operator realization} we have
\begin{align*}
  \pi_{\delta,p}(a)=-\sum_{i=1}^{\ell+1} \bigg[{\ad_{u(x)}e^{\ad_{u(x)}} \over e^{\ad_{u(x)}}-{\rm id}_{\mfrak{n}_+}}\,(e^{-\ad_{u(x)}}a)_{\mfrak{n}_+}\bigg]_i\partial_{u^i} + (\lambda_{\delta,p}+\rho)((e^{-\ad_{u(x)}}a)_{\mfrak{b}})
\end{align*}
for $a \in \mfrak{n}_-$, where $u^i=x_{i-1}$ for $i=1,2,\dots,\ell$ and $u^{\ell+1}=t$, which we can again rewrite using \eqref{eq:Bernoulli polynomials} into the form
\begin{align*}
  \pi_{\delta,p}(a)=-\sum_{i=1}^{\ell+1} \sum_{j=0}^\infty {(-1)^jB_j \over j!}\, [\ad^j_{u(x)}((e^{-\ad_{u(x)}}a)_{\mfrak{n}_+})]_i \partial_{u^i} + (\lambda_{\delta,p} + \rho)((e^{-\ad_{u(x)}}a)_\mfrak{b})
\end{align*}
for $a \in \mfrak{n}_-$, where we used the fact $B_j(1)=(-1)^jB_j(0)$ for $j \in \N_0$.

Hence, we can write
\begin{align*}
  e^{-\ad_{u(x)}}P_n = \sum_{j=0}^\infty {(-1)^j \over j!} \ad^j_{u(x)}(P_n)=\sum_{j=0}^\infty \binom{n}{j}t^j P_{n-j}= \sum_{j=0}^n \binom{n}{j}t^j P_{n-j}= \sum_{j=0}^n \binom{n}{j}t^{n-j}P_j,
\end{align*}
which gives
\begin{align*}
  (e^{-\ad_{u(x)}}P_n)_{\mfrak{n}_+}=\sum_{j=0}^{\ell-1} \binom{n}{j} t^{n-j} P_j \quad \text{and} \quad (e^{-\ad_{u(x)}}P_n)_\mfrak{b}=\sum_{j=\ell}^n \binom{n}{j} t^{n-j} P_j.
\end{align*}
Therefore, we have
\begin{align*}
  \pi_{\delta,p}(P_n)&=-\sum_{i=1}^{\ell+1} \sum_{j=0}^\infty \sum_{k=0}^{\ell-1} {(-1)^jB_j \over j!} \binom{n}{k} t^{n-k} [\ad^j_{u(x)}(P_k)]_i \partial_{u^i}+\sum_{k=\ell}^n \binom{n}{k}t^{n-k}(\lambda_{\delta,p}+\rho)(P_k) \\
  &= -\sum_{j=0}^\infty \sum_{k=0}^{\ell-1} B_j \binom{n}{k} \binom{k}{j} t^{n-k+j} \partial_{x_{k-j}}+  p\binom{n}{\ell}t^{n-\ell} \\
  &= -\sum_{j=0}^{\ell-1} \sum_{k=0}^{\ell-1} B_j \binom{n}{k} \binom{k}{j} t^{n-k+j} \partial_{x_{k-j}}+  p\binom{n}{\ell}t^{n-\ell} \\
  &= -\sum_{j=0}^{\ell-1} \sum_{k=n-\ell+1}^{n-j} B_j \binom{n}{k}\!\binom{n-k}{j} t^{k+j} \partial_{x_{n-k-j}} +  p\binom{n}{\ell}t^{n-\ell}.
\end{align*}
From \eqref{eq:ad(u_e(x)) relation k=1} and \eqref{eq:ad(u_e(x)) relation C,D}, we get
\begin{align*}
  (e^{-\ad_{u(x)}}C)_\mfrak{b}=C+tD+(\ell+1)x_{\ell-1}P_\ell,
\end{align*}
since $(\ad^k_{u(x)}(C))_{\mfrak{n}_+}=0$ for $k>2$, and
\begin{align*}
  (e^{-\ad_{u(x)}}C)_{\mfrak{n}_+}= e^{-\ad_{u(x)}}C- (e^{-\ad_{u(x)}}C)_\mfrak{b} = e^{-\ad_{u(x)}}C- C -tD-(\ell+1)x_{\ell-1}P_\ell.
\end{align*}
Then we may write
\begin{align*}
  -{\ad_{u_e\!(x)}e^{\ad_{u(x)}} \over e^{\ad_{u(x)}}-{\rm id}_{\mfrak{n}_+}}\,(e^{-\ad_{u(x)}}C)_{\mfrak{n}_+}&=  {\ad_{u(x)}e^{\ad_{u(x)}} \over e^{\ad_{u(x)}}-{\rm id}_{\mfrak{n}_+}}\,(C- e^{-\ad_{u(x)}}C + tD +(\ell+1)x_{\ell-1}P_\ell) \\
  &=\ad_{u(x)}(C)+  {\ad_{u(x)}e^{\ad_{u(x)}} \over e^{\ad_{u(x)}}-{\rm id}_{\mfrak{n}_+}}\,(tD +(\ell+1)x_{\ell-1}P_\ell)\\
  &= \ad_{u(x)}(C)+ \sum_{j=0}^\infty {(-1)^jB_j \over j!} \ad^j_{u(x)}(tD+ (\ell+1)x_{\ell-1}P_\ell),
\end{align*}
which allows to simplify the right hand side into
\begin{multline*}
  -\sum_{j=0}^{\ell-2} (2\ell-j)x_jP_{j+1}-\sum_{j=0}^{\ell-1} (\ell-j) tx_j P_j -t^2H + \sum_{j=1}^\infty (\ell+1)B_j \binom{\ell}{j}t^jx_{\ell-1}P_{\ell-j}\\
   + \sum_{j=2}^\infty \sum_{k=0}^{\ell-2} {2(\ell-1-k) \over k+1}B_j \binom{k+1}{j} t^j x_k P_{k-j+1}.
\end{multline*}
This, together with
\begin{align*}
  (\lambda_{\delta,p}+\rho)((e^{-\ad_{u(x)}}C)_\mfrak{b})= \big(\delta-\big(1+{\textstyle {1 \over 2}}\ell(\ell+1)\!\big)\!\big)t+ p(\ell+1)x_{\ell-1}
\end{align*}
results into the final formula
\begin{align*}
  \pi_{\delta,p}(C)&=  -\sum_{j=0}^{\ell-2} (2\ell-j)x_j\partial_{x_{j+1}}-\sum_{j=0}^{\ell-1} (\ell-j) tx_j \partial_{x_j} -t^2\partial_t + \sum_{j=1}^\ell (\ell+1)B_j \binom{\ell}{j}t^jx_{\ell-1} \partial_{x_{\ell-j}}\\
   & \quad + \sum_{j=2}^{\ell-1} \sum_{k=0}^{\ell-2} {2(\ell-1-k)\over k+1}B_j \binom{k+1}{j} t^j x_k \partial_{x_{k-j+1}}
   +  \big(\delta-\big(1+{\textstyle {1 \over 2}}\ell(\ell+1)\!\big)\!\big)t+ p(\ell+1)x_{\ell-1}.
\end{align*}
The computation of the algebraic Fourier transform of all operators is straightforward.}
\vspace{-2mm}


\section{Verma modules for conformal Galilei algebras}
\label{sec:vermmodsingsec}

Based on the construction of their geometric realization, this section presents a detailed analysis of the composition structure of Verma modules $M_\ell(\delta,p)$ for conformal Galilei algebras $\mfrak{cga}_\ell(1,\C)$ with $\ell \in \N$. Namely, we conclude the existence and precise positions of all singular vectors in Verma modules $M_\ell(\delta,p)$ for $\ell \in \N$ and $\delta,p \in \C$. Consequently, we shall describe all homomorphisms
\begin{align}
  \varphi \colon M_\ell(\veps,q) \rarr M_\ell(\delta,p)
\end{align}
of Verma modules for $\veps,q,\delta,p  \in \C$. They are uniquely determined by $\varphi(v_{\veps,q}) \in M_\ell(\delta,p)$, where $v_{\veps,q}$ is the lowest weight vector of $M_\ell(\veps,q)$ with lowest weight $\veps \omega_D+q \omega_{P_\ell}$. By \eqref{eq:Verma module isomorphism}, we have an isomorphism
\begin{align}
  \tau \circ \Phi_{\delta,p} \colon M_\ell(\delta-\delta_\rho,p) \riso \C[\mfrak{n}_+^*]
\end{align}
of $\mfrak{g}$-modules, where the corresponding action of $\mfrak{g}$ on $\C[\mfrak{n}_+^*]$ is given through the mapping
\begin{align}
\hat{\pi}_{\delta,p} \colon \mfrak{g} \rarr \eus{A}^\mfrak{g}_{\mfrak{n}_+^*}.
\end{align}
Now, let $v_{{\rm sing}} \in \C[\mfrak{n}_+^*]$ be a singular vector with a character $\veps \widetilde{\omega}_D+q\widetilde{\omega}_{P_\ell}$. Then from Theorem \ref{thm:cga realization} we obtain that
\begin{align}\label{eq:singular vector condition}
\begin{aligned}
  \hat{\pi}_{\delta,p}(P_\ell)v_{{\rm sing}}=qv_{{\rm sing}}, \quad
  \hat{\pi}_{\delta,p}(D)v_{{\rm sing}}= \veps v_{{\rm sing}},  \quad
  \hat{\pi}_{\delta,p}(C)v_{{\rm sing}}=0, \quad
  \hat{\pi}_{\delta,p}(P_n)v_{{\rm sing}}=0
\end{aligned}
\end{align}
for $n=\ell+1,\ell+2,\dots,2\ell$ and some $\veps,q \in \C$. Since the only eigenvalue of the operator $\hat{\pi}_{\delta,p}(P_\ell)=\ell y_{\ell-1}\partial_u +p$ on $\C[\mfrak{n}_+^*]$ is $p$, we obtain immediately that the singular vector $v_{{\rm sing}} \in \C[\mfrak{n}_+^*]$ does not depend on the variable $u$.

Let us define a Lie subalgebra $\mfrak{n}_{+,0}$ of $\mfrak{n}_+$ and its complement $\mfrak{n}_{+,1}$ in $\mfrak{n}_+$ by
\begin{align}
  \mfrak{n}_{+,0}= \langle P_0,P_1,\dots,P_{\ell-1} \rangle \quad \text{and} \quad \mfrak{n}_{+,1} = \langle H \rangle.
\end{align}
Then the decomposition
\begin{align}
  \mfrak{n}_+ = \mfrak{n}_{+,0} \oplus \mfrak{n}_{+,1}
\end{align}
of the vector space $\mfrak{n}_+$ gives rise to a canonical isomorphism
\begin{align}
  \C[\mfrak{n}^*_+] \simeq \C[\mfrak{n}^*_{+,0}] \otimes_\C \C[\mfrak{n}^*_{+,1}]
\end{align}
of commutative $\C$-algebras. Because $v_{{\rm sing}} \in \C[\mfrak{n}^*_+]$ does not depend on the variable $u$, we have $v_{{\rm sing}} \in \C[\mfrak{n}^*_{+,0}]$ and the system of equations \eqref{eq:singular vector condition} reduces to
\begin{align}
  \hat{\pi}_{\delta,p}^{{\rm res}}(P_\ell)v_{{\rm sing}}=q v_{{\rm sing}}, \quad \hat{\pi}_{\delta,p}^{{\rm res}}(D)v_{{\rm sing}}=\veps v_{{\rm sing}},\quad \hat{\pi}_{\delta,p}^{{\rm res}}(C)v_{{\rm sing}}=0
\end{align}
for some $\veps,q \in \C$, where
\begin{align}
\begin{aligned}
  \hat{\pi}_{\delta,p}^{{\rm res}}(P_\ell)=p, \qquad  \hat{\pi}_{\delta,p}^{{\rm res}}(D)= \sum_{j=0}^{\ell-1} 2(\ell-j)y_j\partial_{y_j} +
    \delta - \delta_\rho, \\
  \hat{\pi}_{\delta,p}^{{\rm res}}(C)=\sum_{j=0}^{\ell-2} (2\ell-j)y_{j+1}\partial_{y_j}-p(\ell+1)\partial_{y_{\ell-1}}.
\end{aligned}
\end{align}
It is more convenient to consider the following change of linear coordinate functions on $\mfrak{n}_{+,0}^*$
defined by
\begin{align}\label{eq:coord change}
 z_n={y_{\ell-1-n} \over (\ell+1+n)!} \quad \text{for} \quad n=0,1,\dots ,\ell-1
\end{align}
with the inverse given by
\begin{align}\label{eq:coord change inv}
 y_n=(2\ell-n)!z_{\ell-1-n} \quad \text{for} \quad n=0,1,\dots ,\ell-1 .
\end{align}
Let us note that the grading of the Lie algebra $\mfrak{g}$ gives us a natural grading on the commutative $\C$-algebra $\C[\mfrak{n}_{+,0}^*]$
defined by $\deg y_n = \ell-n$ for $n=0,1,\dots,\ell-1$. Then \eqref{eq:coord change} implies $\deg z_n=n+1$ for $n=0,1,\dots,\ell-1$.
\medskip

\lemma{\label{soloperz}Let $\delta,p \in \C$. Then as operators in the variables $z_0,z_1,\dots,z_{\ell-1}$, we have
\begin{enumerate}
\item[1)]
$\hat{\pi}^{{\rm res}}_{\delta,p}(P_\ell)=p$,
\item[2)]
$\hat{\pi}^{{\rm res}}_{\delta,p}(D)=\sum\limits_{n=0}^{\ell-1} 2(n+1)z_n\partial_{z_n}+\delta-\delta_\rho$,
\item[3)]
$\hat{\pi}^{{\rm res}}_{\delta,p}(C)=\sum\limits_{n=0}^{\ell-2} z_n\partial_{z_{n+1}} - {p \over \ell!}\partial_{z_0}$.
\end{enumerate}}

\proof{The proof is a straightforward change of coordinates on $\mfrak{n}_{+,0}^*$.}

Let us denote by $F=\C[z_0,z_1,\dots]$ the $\C$-algebra of polynomials in infinitely many variables $z_n$, $n\in \N_0$, and by $F_\ell=\C[z_0,z_1,\dots,z_{\ell-1}]$ for $\ell \in \N$ its finitely generated $\C$-subalgebra of polynomials in the variables $z_n$, $n=0,1,\dots,{\ell-1}$.
Furthermore, let us consider a first order linear differential operator
\begin{align}
  T_1^c=\sum_{n=0}^\infty z_n\partial_{z_{n+1}} - c \partial_{z_0}
\end{align}
on $F$ for $c \in \C$. Then $F_\ell$ for $\ell \in \N$ is clearly an invariant subspace of $T_1^c$ for all $c \in \C$.

Let $v_{{\rm sing}} \in F_\ell \simeq \C[\mfrak{n}_{+,0}^*]$ be a singular vector in $M_\ell(\delta-\delta_\rho,p) \simeq \C[\mfrak{n}_+^*]$ with the character $\veps \widetilde{\omega}_D + q \widetilde{\omega}_{P_\ell}$. Then by Lemma \ref{soloperz} we have
\begin{align}
\hat{\pi}^{{\rm res}}_{\delta,p}(P_\ell)v_{{\rm sing}}=qv_{{\rm sing}} \quad \text{and} \quad \hat{\pi}^{{\rm res}}_{\delta,p}(D)v_{{\rm sing}}=\veps v_{{\rm sing}},
\end{align}
which implies that the singular vector $v_{{\rm sing}}$ is a graded homogeneous polynomial of homogeneity $h={1\over 2}(\veps-\delta+\delta_\rho)$ with respect to the grading $\deg z_n = n+1$ for $n\in \N_0$, and that $q=p$. Finally, the equation $\smash{\hat{\pi}^{{\rm res}}_{\delta,p}}(C)v_{{\rm sing}}=0$ gives $v_{{\rm sing}} \in \ker T_1^c$ with $c = {p \over \ell!}$.

On the other hand, let $v \in F$ be a graded homogeneous polynomial of homogeneity $h \in \N_0$ and let $v \in \ker T_1^c$ for some $c \in \C$. Obviously, there exists a smallest $\ell_0 \in \N$ such that $v \in F_{\ell_0}$. Then for all $\ell \geq \ell_0$, we have $v \in F_\ell$ and
\begin{align}
\begin{gathered}
\hat{\pi}^{{\rm res}}_{\delta,\ell! c}(P_\ell)v=\ell!c\, v, \\
\hat{\pi}_{\delta,\ell!c}^{{\rm res}}(D)v= \bigg({\textstyle \sum\limits_{n=0}^{\ell-1}} 2(n+1)z_n\partial_{z_n} +
    \delta-\delta_\rho\!\bigg)v=(2h+\delta-\delta_\rho)v,\\
  \hat{\pi}^{{\rm res}}_{\delta,\ell!c}(C)v=\bigg({\textstyle \sum\limits_{n=0}^{\ell-2}} z_n\partial_{z_{n+1}} - {\textstyle {\ell!c \over \ell!}}\,\partial_{z_0}\!\bigg)v=0,
\end{gathered}
\end{align}
which means that $v$ is a singular vector in $M_\ell(\delta-\delta_\rho,\ell!c) \simeq \C[\mfrak{n}_+^*]$ for $\mfrak{cga}_\ell(1,\C)$ with character $(2h+\delta-\delta_\rho) \widetilde{\omega}_D + \ell!c\, \widetilde{\omega}_{P_\ell}$.

It follows from our considerations in the last two paragraphs that all singular vectors in the Verma module $M_\ell(\delta-\delta_\rho,p)$ for $\mfrak{cga}_\ell(1,\C)$ with $\ell \in \N$, can be obtained as graded homogeneous solutions of $T_1^c$ for some $c \in \C$ on the vector space $F$. This passage is the key observation with several useful consequences, which we shall use through the rest of the article.


\subsection{Singular vectors and generating series}

The classification of all singular vectors in the Verma module $M_\ell(\delta,p)$ for $\mfrak{cga}_\ell(1,\C)$ with $\ell \in \N$
can be rephrased as a task to determine all solutions of the first order linear differential operator
\begin{align}
  T^c_1=\sum_{n=0}^\infty z_n\partial_{z_{n+1}}-c \partial_{z_0},\quad c\in \C, \label{eq:operator T^c_1}
\end{align}
acting on the vector space $F=\C[z_0,z_1,\dots]$ of polynomials in infinitely many variables $z_n$, $n\in \N_0$. The aim of this section is to describe all solutions of $T^c_1$ on $F$. We remark that the standard notation used in
the rest of the article and related to symmetric polynomials is for the convenience of the reader
reviewed in Appendix \ref{app:symetric polynomials and fucntions}.

We denote by $F^j$ the vector subspace of $F$ of graded homogenous (with respect to the grading $\deg z_n = n+1$ for $n \in \N_0$) polynomials of homogeneity $j \in \N_0$, and by $F^{j,k}$ the vector subspace of $F$ of graded homogenous (with respect to the grading $\deg z_n = n+1$ for $n \in \N_0$) polynomials of homogeneity $j \in \N_0$ and homogeneous (with respect to the grading $\deg z_n =1$ for $n \in \N_0$) polynomials of degree $k \in \N_0$.
Then there are the decompositions
\begin{align}
  F = \bigoplus_{j=0}^\infty F^j \qquad \text{and} \qquad F = \bigoplus_{k=0}^\infty \bigoplus_{j=k}^\infty F^{j,k}.
\end{align}
Now, let us introduce a generating (or formal power) series $a_c(w) \in F[[w]]$ for $c \in \C$ in a formal variable $w$ by
\begin{align}
  a_c(w)=-c + \sum_{n=0}^\infty z_n w^{n+1}.
\end{align}
Further, for $\alpha \in \C$ we have
\begin{align*}
  T^c_1a_c(\alpha w)&= \bigg(\sum_{k=0}^\infty z_k\partial_{z_{k+1}}-c \partial_{z_0}\!\bigg)\!\bigg(\!-c + \sum_{n=0}^\infty z_n (\alpha w)^{n+1}\!\bigg) = -c(\alpha w) + \sum_{k=0}^\infty \sum_{n=0}^\infty z_k \partial_{z_{k+1}}z_n (\alpha w)^{n+1}  \\
  &=  -c(\alpha w) + \sum_{k=0}^\infty \sum_{n=0}^\infty z_k \delta_{n,k+1} (\alpha w)^{n+1} =  -c(\alpha w) + \sum_{k=0}^\infty z_k (\alpha w)^{k+2} = \alpha w a_c(\alpha w),
\end{align*}
which implies
\begin{align}
  T^c_1 \bigg(\prod_{k=1}^r a_c(\alpha_k w)\!\bigg) = \bigg(\sum_{k=1}^r \alpha_k\!\bigg) w \prod_{k=1}^r a_c(\alpha_k w) \label{eq:product series}
\end{align}
for all $r \in \N$ and $\alpha_1,\alpha_2,\dots,\alpha_r \in \C$. Therefore, from \eqref{eq:product series} it follows that the coefficient by $w^n$ in the expansion of $\prod_{k=1}^r\! a_c(\alpha_k w)$ is in $\ker T^c_1$ provided $\sum_{k=1}^r\! \alpha_k = 0$.

Our next step is to find the coefficient by $w^n$ in the expansion of $\prod_{k=1}^r \!a_c(\alpha_k w)$ for $r \in \N$ and $\alpha_1,\alpha_2,\dots,\alpha_r \in \C$. Since this expansion depends on whether or not $c$ is zero, we have to discuss
these two cases separately.
\medskip

\noindent {\bf Case $c \neq 0$.} In order to determine the coefficient by $w^n$, we make use of the following trick. We can formally write
\begin{align}
  -c + \sum_{i=0}^{n-1} z_i w^{i+1} = -c\prod_{i=1}^n (1+\gamma_i w) = -c \sum_{i=0}^n e_i(\gamma_1,\gamma_2,\dots,\gamma_n) w^i,
\end{align}
where
\begin{align}
e_0(\gamma_1,\gamma_2,\dots,\gamma_n)=1,\qquad e_i(\gamma_1,\gamma_2,\dots,\gamma_n)=-{z_{i-1} \over c} \label{eq:expresion e in z c}
\end{align}
for $i=1,2,\dots,n$ are the elementary symmetric polynomials. Therefore, we have
\begin{align}\label{eq:product generating series c}
\begin{aligned}
  \prod_{k=1}^r a_c(\alpha_k w)\ {\rm mod}\ (w^{n+1}) &= \prod_{k=1}^r (-c) \prod_{i=1}^n (1+\gamma_i \alpha_k w)\ {\rm mod}\ (w^{n+1}) \\
  & = (-c)^r \prod_{i=1}^n \prod_{k=1}^r (1 + \gamma_i \alpha_k w)\ {\rm mod}\ (w^{n+1})
\end{aligned}
\end{align}
for all $n \in \N_0$, where we always regard the empty product equal to $1$. By \eqref{eq:Cauchy kernel}, we get
\begin{align}
\begin{aligned}
  \prod_{i=1}^n \prod_{k=1}^r (1 + \gamma_i \alpha_k w)  = \sum_{\lambda \in \mcal{P}} \veps_\lambda {p_\lambda(\alpha_1,\dots,\alpha_r) p_\lambda(\gamma_1,\dots,\gamma_n) \over z_\lambda}\, w^{|\lambda|},
\end{aligned}
\end{align}
which enables us to rewrite \eqref{eq:product generating series c} into the form
\begin{align}
  \prod_{k=1}^r a_c(\alpha_k w)\ {\rm mod}\ (w^{n+1}) = (-c)^r \sum_{k=0}^n \sum_{\lambda \in \mcal{P}_k} \veps_\lambda {p_\lambda(\alpha_1,\dots,\alpha_r) p_\lambda(\gamma_1,\dots,\gamma_n) \over z_\lambda}\, w^k.
\end{align}
Using \eqref{eq:expresion e in z c} and the Newton's identities \eqref{eq:expression p in e}, the polynomial $p_k(\gamma_1,\dots,\gamma_n)$ can be written as a graded homogenous (with respect to the grading $\deg z_n = n+1$ for $n \in \N_0$) polynomial of homogeneity $k$ in the variables $z_0, z_1,\dots,z_{k-1}$ provided $k \leq n$. Hence, we have
\begin{align}
  p_k(\gamma_1,\dots,\gamma_n) = u_k\Big(\!-\!{z_0 \over c},\dots,-{z_{k-1}\over c}\Big) = {t^c_k(z) \over (-1)^kc^k}  \label{eq:t^c_k(z) definition}
\end{align}
as follows from \eqref{eq:relation p and e}, where the polynomial $t^c_k(z) \in F^k$ and does not depend on $n$. Let us note that the polynomial $t^c_k(z)$ has a well-defined limit for $c \rarr 0$, and we get $t^0_k(z)=z_0^k$ for $k \in \N$. Therefore, we may write
\begin{align}
  \prod_{k=1}^r a_c(\alpha_k w)\ {\rm mod}\ (w^{n+1}) = (-c)^r\sum_{\lambda \in \mcal{P}} \veps_\lambda {p_\lambda(\alpha_1,\dots,\alpha_r) t^c_\lambda(z) \over z_\lambda}\, {w^{|\lambda|} \over (-1)^{|\lambda|}c^{|\lambda|}}\,\ {\rm mod}\ (w^{n+1})
\end{align}
for all $n\in \N_0$, where $t^c_\lambda(z)=\prod_{i \in \N} t^c_i(z)^{m_i(\lambda)}$, which gives
\begin{align} \label{eq:expansion product in t c}
\begin{aligned}
  \prod_{k=1}^r a_c(\alpha_k w) = (-c)^r\sum_{\lambda \in \mcal{P}} \veps_\lambda {p_\lambda(\alpha_1,\dots,\alpha_r) t^c_\lambda(z) \over z_\lambda}\, {w^{|\lambda|} \over (-1)^{|\lambda|}c^{|\lambda|}}.
\end{aligned}
\end{align}
We obtain the following power series expansion
\begin{align}
  \prod_{k=1}^r a_c(\alpha_k w) = \sum_{n=0}^\infty \sum_{\lambda \in \mcal{P}_n} p_\lambda(\alpha_1,\dots,\alpha_r)  {\veps_\lambda t^c_\lambda(z) \over z_\lambda}\, {w^n \over (-1)^{n+r} c^{n+r}} \label{eq:expansion c}
\end{align}
for all $r \in \N$ and $\alpha_1,\alpha_2,\dots,\alpha_r \in \C$.
\medskip

\noindent {\bf Case $c = 0$.} We can again formally write
\begin{align}
  \sum_{i=0}^n z_i w^{i+1} = z_0w \prod_{i=1}^n (1+\beta_i  w)= z_0w \sum_{i=0}^n e_i(\beta_1,\beta_2,\dots,\beta_n) w^i,
\end{align}
where
\begin{align}
e_i(\beta_1,\beta_2,\dots,\beta_n)={z_i \over z_0} \label{eq:expresion e in z}
\end{align}
for $i=0,1,\dots,n$ are the elementary symmetric polynomials. Therefore, we have
\begin{align}\label{eq:product generating series}
\begin{aligned}
  \prod_{k=1}^r a_0(\alpha_k w)\ {\rm mod}\ (w^{n+2}) &= \prod_{k=1}^r z_0w \prod_{i=1}^n (1+\beta_i \alpha_k w)\ {\rm mod}\ (w^{n+2}) \\
  & = (z_0w)^r \prod_{i=1}^n \prod_{k=1}^r (1 + \beta_i \alpha_k w)\ {\rm mod}\ (w^{n+2})
\end{aligned}
\end{align}
for all $n \in \N_0$. Now, from \eqref{eq:Cauchy kernel} we get
\begin{align}
\begin{aligned}
  \prod_{i=1}^n \prod_{k=1}^r (1 + \beta_i \alpha_k w) = \sum_{\lambda \in \mcal{P}} \veps_\lambda {p_\lambda(\alpha_1,\dots,\alpha_r) p_\lambda(\beta_1,\dots,\beta_n) \over z_\lambda}\, w^{|\lambda|},
\end{aligned}
\end{align}
which enables us to rewrite \eqref{eq:product generating series} into the form
\begin{align}
  \prod_{k=1}^r a_0(\alpha_k w)\  {\rm mod}\ (w^{n+1+r}) = (z_0w)^r \sum_{k=0}^n \sum_{\lambda \in \mcal{P}_k} \veps_\lambda {p_\lambda(\alpha_1,\dots,\alpha_r) p_\lambda(\beta_1,\dots,\beta_n) \over z_\lambda}\, w^k.
\end{align}
Using \eqref{eq:expresion e in z} and the Newton's identities \eqref{eq:expression p in e}, the polynomial $p_k(\beta_1,\dots,\beta_n)$ can be written as a graded homogenous (with respect to the grading $\deg {z_n \over z_0}= n$ for $n \in \N$) polynomial of homogeneity $k$ in the variables ${z_1 \over z_0},{z_2 \over z_0},\dots,{z_k \over z_0}$ provided $k \leq n$. Hence, we have
\begin{align}
  p_k(\beta_1,\dots,\beta_n) = u_k\Big({z_1 \over z_0},\dots,{z_k \over z_0}\Big) = {t_k(z) \over z_0^k} \label{eq:t_k(z) definition}
\end{align}
as follows from \eqref{eq:relation p and e}, where the polynomial $t_k(z) \in F^{2k,k}$ and does not depend on $n$. Therefore, we may write
\begin{align}
  \prod_{k=1}^r a_0(\alpha_k w)\ {\rm mod}\ (w^{n+1+r}) = (z_0w)^r\sum_{\lambda \in \mcal{P}} \veps_\lambda {p_\lambda(\alpha_1,\dots,\alpha_r) t_\lambda(z) \over z_\lambda}\, {w^{|\lambda|} \over z_0^{|\lambda|}}\,\ {\rm mod}\ (w^{n+1+r})
\end{align}
for all $n\in \N_0$, where $t_\lambda(z)=\prod_{i \in \N} t_i(z)^{m_i(\lambda)}$, which gives
\begin{align} \label{eq:expansion product in t}
\begin{aligned}
  \prod_{k=1}^r a_0(\alpha_k w) = (z_0w)^r\sum_{\lambda \in \mcal{P}} \veps_\lambda {p_\lambda(\alpha_1,\dots,\alpha_r) t_\lambda(z) \over z_\lambda}\, {w^{|\lambda|} \over z_0^{|\lambda|}}.
\end{aligned}
\end{align}
Further, for $\lambda \in \mcal{P}_n$ we have
\begin{align}
  p_\lambda(\alpha_1,\dots,\alpha_r)= \sum_{\mu \in \mcal{P}_n\!(r)} a^r_{\lambda,\mu} p_\mu(\alpha_1,\dots,\alpha_r),
\end{align}
where $a^r_{\lambda,\mu} \in \C$. Therefore, from \eqref{eq:expansion product in t} we get
\begin{align*}
  \prod_{k=1}^r a_0(\alpha_k w) &= (z_0w)^r\sum_{n=0}^\infty \sum_{\lambda \in \mcal{P}_n} \veps_\lambda {p_\lambda(\alpha_1,\dots,\alpha_r) t_\lambda(z) \over z_\lambda}\, {w^n \over z_0^n} \\
  &= (z_0w)^r\sum_{n=0}^\infty \sum_{\lambda \in \mcal{P}_n} \sum_{\mu \in \mcal{P}_n\!(r)} \veps_\lambda {a^r_{\lambda,\mu}p_\mu(\alpha_1,\dots,\alpha_r) t_\lambda(z) \over z_\lambda}\, {w^n \over z_0^n} \\
  &= \sum_{n=0}^\infty \sum_{\mu \in \mcal{P}_n\!(r)}  p_\mu(\alpha_1,\dots,\alpha_r) \bigg(\sum_{\lambda \in \mcal{P}_n}   {\veps_\lambda a^r_{\lambda,\mu}  \over z_\lambda}\,t_\lambda(z) {z_0^r \over z_0^n} \bigg) w^{n+r},
\end{align*}
and denoting
\begin{align}
  s^r_\lambda(z)= \sum_{\mu \in \mcal{P}_n}   {\veps_\mu a^r_{\mu,\lambda}  \over z_\mu}\,t_\mu(z) {z_0^r \over z_0^n} \label{eq:s_lambda formula}
\end{align}
for $\lambda \in \mcal{P}_n(r)$, we obtain the following power series expansion
\begin{align}
  \prod_{k=1}^r a_0(\alpha_k w) = \sum_{n=0}^\infty \sum_{\lambda \in \mcal{P}_n\!(r)} p_\lambda(\alpha_1,\dots,\alpha_r) s^r_\lambda(z)  w^{n+r} \label{eq:expansion}
\end{align}
for all $r \in \N$ and $\alpha_1,\alpha_2,\dots,\alpha_r \in \C$.
\medskip

\proposition{\label{prop:linear independence}
\begin{enum}
  \item[1)] If $c \neq 0$, then the set of polynomials $\{t^c_\lambda(z);\, \lambda \in \mcal{P}\}$  in $F$ is linearly independent.
  \item[2)] Let us define $s^0_\emptyset(z)=1$. Then the set of polynomials $\{s^r_\lambda(z);\, r \in \N_0,\, \lambda \in \mcal{P}(r)\}$  in $F$ is linearly independent.
\end{enum}}

\proof{1) We have the decomposition \smash{$F = \bigoplus_{j=0}^\infty F^j$}. Since $t^c_\lambda(z) \in F^{|\lambda|}$, which follows from the fact $t^c_\lambda(z)= \smash{\prod_{i \in \N} t^c_i(z)^{m_i(\lambda)}}$ and $t^c_k(z) \in F^k$, it is enough to show that the set of polynomials $\{t^c_\lambda(z);\, \lambda \in \mcal{P}_n\}$  is linearly independent for all $n \in \N_0$. Let us suppose that \smash{$\sum_{\lambda \in \mcal{P}_n}\!a_\lambda t^c_\lambda(z)=0$} for some $a_\lambda \in \C$. If $c \neq 0$, then using \eqref{eq:t^c_k(z) definition}, we may write
\begin{align*}
  0= \sum_{\lambda \in \mcal{P}_n} a_\lambda t^c_\lambda(z) = \sum_{\lambda \in \mcal{P}_n} (-1)^{|\lambda|}c^{|\lambda|} a_\lambda\, {t^c_\lambda(z) \over (-1)^{|\lambda|}c^{|\lambda|}} = \sum_{\lambda \in \mcal{P}_n} (-1)^{|\lambda|}c^{|\lambda|} a_\lambda p_\lambda(\gamma_1,\gamma_2,\dots,\gamma_n),
\end{align*}
where
\begin{align*}
e_0(\gamma_1,\gamma_2,\dots,\gamma_n)=1,\qquad e_i(\gamma_1,\gamma_2,\dots,\gamma_n)=-{z_{i-1} \over c}
\end{align*}
for $i=1,2,\dots,n$. But the set of polynomials $\{p_\lambda(\gamma_1,\gamma_2,\dots,\gamma_n);\, \lambda \in \mcal{P}_n\}$ is linearly independent, hence we have $a_\lambda=0$ for all $\lambda \in \mcal{P}_n$, and we are done.

2) We have the decomposition \smash{$F = \bigoplus_{k=0}^\infty \bigoplus_{j=k}^\infty F^{j,k}$}. Since we have $s^r_\lambda(z) \in \smash{F^{r+|\lambda|,r}}$ for $r \in \N_0$ and $\lambda \in \mcal{P}(r)$, which follows by \eqref{eq:s_lambda formula} and from the fact $t_\lambda(z)= \prod_{i \in \N} \smash{t_i(z)^{m_i(\lambda)}}$ and $t_k(z) \in F^{2k,k}$, it is enough to show that the set of polynomials $\{s^r_\lambda(z);\, \lambda \in \mcal{P}_n(r)\}$  is linearly independent for all $r,n\in \N_0$. As $a^r_{\mu,\lambda}=\delta_{\mu,\lambda}$ for $\mu,\lambda \in \mcal{P}_n(r)$, we obtain
\begin{align*}
  s^r_\lambda(z)= {\veps_\lambda \over z_\lambda}t_\lambda(z) {z_0^r \over z_0^n} + \sum_{\mu \in \mcal{P}_n \setminus \mcal{P}_n\!(r)}   {\veps_\mu a^r_{\mu,\lambda}  \over z_\mu}\,t_\mu(z) {z_0^r \over z_0^n}
\end{align*}
for $\lambda \in \mcal{P}_n(r)$. But the set $\{t_\mu(z);\, \mu \in \mcal{P}_n\}$ is linearly independent by a similar argument as in the first part, therefore the set $\{s^r_\lambda(z);\, \lambda \in \mcal{P}_n(r)\}$ is also linearly independent, and we are done.}

Now, if $\sum_{k=1}^r \!\alpha_k =0$, then the power series expansions \eqref{eq:expansion c} and \eqref{eq:expansion} reduce into
\begin{align}
  \prod_{k=1}^r a_c(\alpha_k w) =
  \begin{cases}
  \sum\limits_{n=0}^\infty \sum\limits_{\substack{\lambda \in \mcal{P}_n \\ m_1(\lambda)=0}} p_\lambda(\alpha_1,\dots,\alpha_r)  {\veps_\lambda t^c_\lambda(z) \over z_\lambda}\, {w^n \over (-1)^{n+r} c^{n+r}} & \text{for $c \neq 0$}, \\[4mm]
  \sum\limits_{n=0}^\infty \sum\limits_{\substack{\lambda \in \mcal{P}_n\!(r) \\ m_1(\lambda)=0}} p_\lambda(\alpha_1,\dots,\alpha_r) s^r_\lambda(z) w^{n+r} & \text{for $c=0$},
  \end{cases}
\end{align}
and as a consequence of \eqref{eq:product series} we obtain that
\begin{align} \label{eq:kernel T1}
   \sum_{\substack{\lambda \in \mcal{P}_n \\ m_1(\lambda)=0}} p_\lambda(\alpha_1,\dots,\alpha_r)  {\veps_\lambda t^c_\lambda(z) \over z_\lambda} \in \ker T_1^{c\neq 0} \quad \text{and} \quad
   \sum_{\substack{\lambda \in \mcal{P}_n\!(r)\\ m_1(\lambda)=0}} p_\lambda(\alpha_1,\dots,\alpha_r) s^r_\lambda(z) \in \ker T^{c=0}_1
\end{align}
for all $n \in \N_0$, $r \in \N$ and $\alpha_1,\alpha_2,\dots,\alpha_r \in \C$ satisfying $\sum_{k=1}^r\! \alpha_k=0$. In fact, we shall prove that $s_\lambda^r(z) \in \ker T^0_1$ for all $r \in \N$ and $\lambda \in \mcal{P}(r)$ satisfying $m_1(\lambda)=0$, and that $t^c_\lambda(z) \in \ker T^c_1$ for all $\lambda \in \mcal{P}$ satisfying $m_1(\lambda)=0$ provided $c \neq 0$.
\medskip

\proposition{\label{prop:power sum system}The following system of equations
\begin{align} \label{eq:power sum system}
\begin{aligned}
 p_1(\alpha_1,\alpha_2,\dots,\alpha_r)&=a_1, \\
 p_2(\alpha_1,\alpha_2,\dots,\alpha_r)&=a_2, \\
 \vdots \\
 p_r(\alpha_1,\alpha_2,\dots,\alpha_r)&=a_r \\
\end{aligned}
\end{align}
has a unique solution up to a permutation of $\alpha_1,\alpha_2,\dots,\alpha_r$ for any $a_1,a_2,\dots,a_r \in \C$.}

\proof{We see that the system of equations \eqref{eq:power sum system} for $a_1,a_2,\dots,a_r \in \C$  is equivalent to the system of equations
\begin{align*}
\begin{aligned}
 e_1(\alpha_1,\alpha_2,\dots,\alpha_r)&=b_1, \\
 e_2(\alpha_1,\alpha_2,\dots,\alpha_r)&=b_2, \\
 \vdots \\
 e_r(\alpha_1,\alpha_2,\dots,\alpha_r)&=b_r \\
\end{aligned}
\end{align*}
for $b_1,b_2,\dots,b_r \in \C$, where $a_1,a_2,\dots,a_r$ and $b_1,b_2,\dots,b_r$ are related through the formulas \eqref{eq:expression p in e} and \eqref{eq:expression e in p}.

Let us consider a polynomial $p(x)=\sum_{i=0}^r (-1)^ib_ix^{r-i}$ with $b_0=1$. Let $\alpha_1,\alpha_2,\dots,\alpha_r \in \C$ be the roots of the polynomial $p(x)$, then we have
\begin{align*}
  \sum_{i=0}^r (-1)^ib_ix^{r-i}= \prod_{i=1}^r (x-\alpha_i)= \sum_{i=0}^r (-1)^ie_i(\alpha_1,\alpha_2,\dots,\alpha_n)x^{r-i},
\end{align*}
which proves the claim.}

\lemma{\label{lem:linear combination}Let us assume that $f_0(z),f_1(z),\dots,f_n(z) \in F$. If $\sum_{i=0}^n \!a^i f_i(z) \in \ker T^c_1$ for all $a \in \C$, then $f_i(z) \in \ker T^c_1$ for $i=0,1,\dots,n$.}

\proof{We introduce $g_a(z)= \sum_{i=0}^n\! a^i f_i(z)$ for $a \in \C$. Then we have
\begin{align*}
  \begin{pmatrix}
    1 & a_0 & a_0^2 & \dots & a_0^n \\
    1 & a_1 & a_1^2 & \dots & a_1^n \\
    \vdots & \vdots & \vdots & \ddots & \vdots \\
    1 & a_n & a_n^2 & \dots & a_n^n \\
  \end{pmatrix}\!
  \begin{pmatrix}
    f_0(z) \\
    f_1(z) \\
    \vdots \\
    f_n(z)
  \end{pmatrix}=
  \begin{pmatrix}
    g_{a_0}(z) \\
    g_{a_1}(z) \\
    \vdots \\
    g_{a_n}(z)
  \end{pmatrix}
\end{align*}
for $a_0,a_1,\dots,a_n \in \C$. If $a_i \neq a_j$ for $i\neq j$, then the Vandermonde matrix is invertible. Therefore, we can express $f_i(z)$ as a $\C$-linear combination of $g_{a_0}(z),g_{a_1}(z),\dots,g_{a_n}(z)$ for $i=0,1,\dots,n$. Because $g_{a_i}(z) \in \ker T^c_1$ for $i=0,1,\dots,n$, we obtain that $f_i(z) \in \ker T^c_1$ for $i=0,1,\dots,n$.}

Let us denote $a^\lambda= \prod_{i=1}^r\! a_i^{m_i(\lambda)}$ for $a_1,a_2,\dots,a_r \in \C$ and $\lambda \in \mcal{P}(r)$. Now, if $c=0$, then by Proposition \ref{prop:power sum system} and \eqref{eq:kernel T1} we obtain that
\begin{align}
  \sum_{\substack{\lambda \in \mcal{P}_n\!(r) \\ m_1(\lambda)=0}} a^\lambda s^r_\lambda(z) \in \ker T^0_1
\end{align}
for all $a_1,a_2,\dots,a_r \in \C$. Since we have
\begin{align*}
\sum_{\substack{\lambda \in \mcal{P}_n\!(r) \\ m_1(\lambda)=0}} a^\lambda s^r_\lambda(z)=  \sum_{k_2=0}^n a_2^{k_2}\bigg( \sum_{k_3=0}^n a_3^{k_3}\bigg( \dots \bigg(\sum_{k_r=0}^n a_r^{k_r} f^r_{0,k_2,\dots,k_r}(z)\!\bigg)\dots\bigg)\!\bigg) \in \ker T^0_1,
\end{align*}
where $\smash{f^r_{k_1,k_2,\dots,k_r}}(z)=s_\lambda^r(z)$ for a partition $\lambda \in \mcal{P}_n(r)$ satisfying $m_i(\lambda)=k_i$ for $i=1,2,\dots,r$, and $f^r_{k_1,k_2,\dots,k_r}(z)=0$ otherwise, the multiple application of Lemma \ref{lem:linear combination} gives $s_\lambda^r(z) \in \ker T^0_1$ for all $\lambda \in \mcal{P}_n(r)$ satisfying $m_1(\lambda)=0$. Further, if $c\neq 0$, then by Proposition \ref{prop:power sum system} and \eqref{eq:kernel T1} we obtain that
\begin{align}
 \sum_{\substack{\lambda \in \mcal{P}_n \\ m_1(\lambda)=0}} a^\lambda  {\veps_\lambda t^c_\lambda(z) \over z_\lambda} \in \ker T_1^c
\end{align}
for all $a_1,a_2,\dots,a_r \in \C$ provided $r \geq n$. By the same argument as for $c =0$, we get $t^c_\lambda(z) \in \ker T^c_1$ for all $\lambda \in \mcal{P}_n$ satisfying $m_1(\lambda)=0$.
\medskip

The last step is to show that we have found all linearly independent solutions of $T^c_1$ on $F$. To prove this fact, let us consider three first order linear differential operators
\begin{align}
  T^c_1&=\sum_{n=0}^\infty z_n\partial_{z_{n+1}}-c \partial_{z_0}, \\
  T^c_0 & = \sum_{n=0}^\infty (n+1) z_n \partial_{z_n} + {1\over 2} c\widebar{c}, \\
  T^c_{-1} & = \sum_{n=0}^\infty (n+1)(n+2) z_{n+1} \partial_{z_n} - \widebar{c}z_0
\end{align}
regarded as endomorphisms of the vector space $F$.
\medskip

\lemma{\label{lem:sl(2,C)}We have the following nontrivial commutation relations
\begin{align}
  [T^c_0, T^c_1]=-T^c_1, \qquad [T^c_1,T^c_{-1}] = 2T^c_0, \qquad [T^c_0,T^c_{-1}]= T^c_{-1}.
\end{align}}

\proof{A straightforward computation gives
\begin{align*}
  [T^c_1,T^c_{-1}] & = \sum_{n,k \in \N_0} (k+1)(k+2)[z_n\partial_{z_{n+1}}, z_{k+1}\partial_{z_k}]+c\widebar{c}\\
   & =  \sum_{n,k \in \N_0} (k+1)(k+2)\big([z_n\partial_{z_{n+1}}, z_{k+1}]\partial_{z_k}+z_{k+1} [z_n\partial_{z_{n+1}},\partial_{z_k}]\big)+ c\widebar{c} \\
   &= \sum_{n,k \in \N_0} (k+1)(k+2)(\delta_{n,k} z_n\partial_{z_k}-\delta_{n,k} z_{k+1} \partial_{z_{n+1}})+c\widebar{c} \\
   &= \sum_{n \in \N_0} ((n+1)(n+2)-n(n+1))z_n\partial_{z_n}+c\widebar{c} = \sum_{n \in \N_0}\!2(n+1)z_n\partial_{z_n}+c\widebar{c}= 2T^c_0.
\intertext{Similarly, we have}
  [T^c_0,T^c_{-1}]& = \sum_{n,k \in \N_0} (k+1)(n+1)(n+2)[z_k \partial_{z_k},z_{n+1} \partial_{z_n}]-\widebar{c}z_0 \\
  & =  \sum_{n,k \in \N_0} (k+1)(n+1)(n+2)\big([z_k \partial_{z_k},z_{n+1}]\partial_{z_n}+ z_{n+1}[z_k \partial_{z_k},\partial_{z_n}] \big) -\widebar{c}z_0\\
  &=\sum_{n,k \in \N_0} (k+1)(n+1)(n+2)(\delta_{n+1,k}z_k\partial_{z_n}- \delta_{n,k}z_{n+1}\partial_{z_k})-\widebar{c}z_0 \\
  &=\sum_{n \in N_0} (n+1)(n+2)((n+2)-(n+1))z_{n+1}\partial_{z_n}-\widebar{c}z_0 \\
  &= \sum_{n \in \N_0} (n+1)(n+2)z_n\partial_{z_{n+1}}-\widebar{c}z_0=T^c_{-1}
\intertext{and}
  [T^c_0,T^c_1]& = \sum_{n,k \in \N_0} (k+1)[z_k \partial_{z_k},z_n \partial_{z_{n+1}}]+c\partial_{z_0}\\
  &=  \sum_{n,k \in \N_0} (k+1)\big([z_k \partial_{z_k},z_n]\partial_{z_{n+1}}+ z_n[z_k \partial_{z_k},\partial_{z_{n+1}}] \big) +c\partial_{z_0}\\
  &=\sum_{n,k \in \N_0} (k+1)(\delta_{n,k}z_k\partial_{z_{n+1}}- \delta_{n+1,k}z_n\partial_{z_k})+c\partial_{z_0}\\
  &=\sum_{n \in N_0} ((n+1)-(n+2))z_n\partial_{z_{n+1}}+c\partial_{z_0} = -\sum_{n \in \N_0} z_n\partial_{z_{n+1}}+c\partial_{z_0}=-T^c_1.
\end{align*}
This finishes the proof.}

Now, from Lemma \ref{lem:sl(2,C)} it follows that the Lie algebra generated by the linear differential operators $T^c_{-1}$, $T^c_0$ and $T^c_1$ is isomorphic to the Lie algebra $\mfrak{sl}(2,\C)$. Let us denote by
\begin{align}
  e=\begin{pmatrix}
    0 & 1 \\
    0 & 0
  \end{pmatrix}\!, \qquad
  h=\begin{pmatrix}
    1 &  0\\
    0 & -1
  \end{pmatrix}\!, \qquad
  f=\begin{pmatrix}
    0 &  0\\
    1 & 0
  \end{pmatrix}
\end{align}
the standard basis of $\mfrak{sl}(2,\C)$. Then the subspace $\mfrak{h}=\C h$ is a Cartan subalgebra and the subspace $\mfrak{b}=\C h \oplus \C e$ is a Borel subalgebra of $\mfrak{sl}(2,\C)$. Let us define $\alpha \in \mfrak{h}^*$ by $\alpha(h)=2$ and $\omega \in \mfrak{h}^*$ by $\omega(h)=1$. The representation $\sigma_c \colon \mfrak{sl}(2,\C) \rarr \mfrak{gl}(F)$ of $\mfrak{sl}(2,\C)$ is given by
\begin{align}
  \sigma_c(e)=T^c_1,\qquad \sigma_c(h)=-2T^c_0, \qquad \sigma_c(f)=-T^c_{-1}.
\end{align}
Hence, we have got an infinite-dimensional representation of $\mfrak{sl}(2,\C)$ on the vector space $F$. Moreover, if we introduce an inner product on $F$ by
\begin{align}\label{innprod}
  \langle g_1, g_2 \rangle = g^*_1(\tilde{\partial}_z)g_2(z)|_{z=0},
\end{align}
where $z=(z_0,z_1,\dots,z_n,\dots)$, $\tilde{\partial}_z=(\partial_{z_0},{1 \over 2}\partial_{z_1},\dots,{1 \over n!(n+1)!}\partial_{z_n},\dots$) and $g^*_1(z)=\overline{g_1(\widebar{z})}$, then we have
\begin{align}
  \langle T^c_k g_1, g_2 \rangle = \langle g_1, T^c_{-k} g_2 \rangle
\end{align}
for all $g_1, g_2 \in F$ and $k=-1,0,1$. The existence of a contravariant inner product implies that the representations $\sigma_c \colon \mfrak{sl}(2,\C) \rarr \mfrak{gl}(F)$ of $\mfrak{sl}(2,\C)$ are completely reducible for all $c\in \C$.

Furthermore, for $\lambda \in \C$ we denote by $M_\mfrak{b}(\lambda\omega)$ the Verma $\mfrak{sl}(2,\C)$-module with highest weight $\lambda \omega \in \mfrak{h}^*$, and by $L_\mfrak{b}(\lambda \omega)$ its quotient by the maximal $\mfrak{sl}(2,\C)$-submodule. As a consequence, we obtain that the $\mfrak{sl}(2,\C)$-submodule of $F$ generated by a polynomial $g \in F$ satisfying $\sigma_c(h)g=\lambda g$ with $\lambda \in \C$ and $\sigma_c(e)g=0$ is isomorphic to $L_\mfrak{b}(\lambda \omega)$.
\medskip

\noindent{\bf Case $c \neq 0$.} Applying $\sigma_c(h)$ to $t^c_\lambda(z)$ for $\lambda \in \mcal{P}$, we obtain
\begin{align}
  \sigma_c(h)t^c_\lambda(z)=-(2|\lambda|+|c|^2)t^c_\lambda(z), \label{eq:t^c_lambda(z) weight}
\end{align}
since $t^c_\lambda(z)= \prod_{i \in \N} t^c_i(z)^{m_i(\lambda)}$ and $t^c_k(z) \in F^k$. Moreover, if $\lambda \in \mcal{P}$ satisfies $m_1(\lambda)=0$, then we also have
\begin{align}
  \sigma_c(e)t^c_\lambda(z)=0.
\end{align}
Therefore, the polynomial $t^c_\lambda(z)$ is the highest weight vector of highest weight $-(2|\lambda|+|c|^2)$ for all $\lambda \in \mcal{P}$ satisfying $m_1(\lambda)=0$.

Let us denote by $m(j)$ for $j\in \N_0$ the dimension of the linear span of $\{t^c_\lambda(z);\, \lambda \in \mcal{P},\, m_1(\lambda)=0,\, |\lambda|=j\}$. Since $t^c_\lambda(z)$ for $\lambda \in \mcal{P}$ are linearly independent as follows from Proposition \ref{prop:linear independence}, we obtain that
\begin{align}
  m(j) = \# \{\lambda \in \mcal{P};\, m_1(\lambda)=0,\, |\lambda|=j\}
\end{align}
for $j \in \N_0$. Hence, we have
\begin{align}
  m(j)= \#\{(m_2,m_3,\dots,m_j)\in \N_0^{j-1};\, {\textstyle \sum_{i=2}^j} im_i=j\}
\end{align}
for $j \in \N_0$.

Let $M$ the $\mfrak{sl}(2,\C)$-submodule of $F$ generated by $\{t^c_\lambda(z);\, \lambda\in \mcal{P},\, m_1(\lambda)=0\}$. We take an orthogonal basis of the linear span of $\{t^c_\lambda(z);\, \lambda\in \mcal{P},\, m_1(\lambda)=0\}$ consisting of weight vectors, so that each basis vector generates an $\mfrak{sl}(2,\C)$-submodule isomorphic to $L_\mfrak{b}(\lambda \omega)$ for some $\lambda \in \C$. Moreover, these $\mfrak{sl}(2,\C)$-submodules are mutually orthogonal for the inner
product \eqref{innprod}, hence we get
\begin{align}
   F \supset M \simeq  M_\mfrak{b}(-|c|^2\omega) \oplus \bigoplus_{j=2}^\infty \bigoplus_{\substack{r_2,r_3,\dots \in \N_0 \\ 2r_2+3r_3+ \dots = j}} M_b(-(2j+|c|^2)\omega), \label{eq:submodule decomposition c not 0}
\end{align}
where the direct summands are generated by $t^c_\emptyset(z)$ and $\{t^c_\lambda(z);\, \lambda \in \mcal{P},\, m_1(\lambda)=0,\, |\lambda|=j\}$ for $j\geq 2$.
\medskip

\noindent{\bf Case $c=0$.}
Applying $\sigma_0(h)$ to $s^r_\lambda(z)$ for $r \in \N_0$ and $\lambda \in \mcal{P}(r)$, we obtain from \eqref{eq:s_lambda formula} that
\begin{align}
  \sigma_0(h)s^r_\lambda(z)= -2(r+|\lambda|)s^r_\lambda(z), \label{eq:s^r_lambda(z) weight}
\end{align}
since $t_\lambda(z)= \prod_{i \in \N} t_i(z)^{m_i(\lambda)}$ and $t_k(z) \in F^{2k,k}$. Moreover, if $\lambda \in \mcal{P}(r)$ satisfies $m_1(\lambda)=0$, then we also have
\begin{align}
  \sigma_0(e)s^r_\lambda(z)=0.
\end{align}
Therefore, the polynomial $s^r_\lambda(z)$ is the highest weight vector of highest weight $-2(r+|\lambda|)$ for all $r \in \N_0$ and $\lambda \in \mcal{P}(r)$ satisfying $m_1(\lambda)=0$.

Let us denote by $m(j)$ for $j\in \N_0$ the dimension of the linear span of $\{s^r_\lambda(z);\, r\in \N_0,\, \lambda \in \mcal{P}(r),\, m_1(\lambda)=0,\, r+|\lambda|=j\}$. Since $s^r_\lambda(z)$ for $r \in \N_0$, $\lambda \in \mcal{P}(r)$ are linearly independent as follows from Proposition \ref{prop:linear independence}, we obtain that
\begin{align}
  m(j) = \# \{(r,\lambda);\, r\in \N_0,\,\lambda \in \mcal{P}(r),\, m_1(\lambda)=0,\, r+|\lambda|=j\}
\end{align}
for $j\in \N_0$. Further, for $j,r \in \N_0$ let us denote
\begin{align}
  m(j,r)=\# \{\lambda \in \mcal{P}(r);\, m_1(\lambda)=0,\, r+|\lambda|=j\}.
\end{align}
Then we have $m(0)=\#\{(0,\emptyset)\}=1$ and $m(1)=\#\{(1,\emptyset)\}=1$. If $j \geq 2$, then $m(j,0)=0$ and $m(j,1)=0$, which gives
\begin{align}
  m(j)=\sum_{r=2}^j m(j,r).
\end{align}
For $j \geq r \geq 2$ we may write
\begin{align*}
  m(j,r)&= \#\{(m_2,m_3,\dots,m_j)\in \N_0^{j-1};\, {\textstyle \sum_{i=2}^r} im_i=j-r,\, m_i=0\ \text{for}\ i>r\} \\
    &= \#\{(m_2,m_3,\dots,m_j)\in \N_0^{j-1};\, {\textstyle \sum_{i=2}^r} im_i=j,\, m_r \geq 1,\, m_i=0\ \text{for}\ i>r\},
\end{align*}
which implies
\begin{align}
  m(j)=\#\{(m_2,m_3,\dots,m_j)\in \N_0^{j-1};\, {\textstyle \sum_{i=2}^j} im_i=j\}
\end{align}
for $j \geq 2$.

Let $M$ be the $\mfrak{sl}(2,\C)$-submodule of $F$ generated by $\{s^r_\lambda(z);\, r \in \N_0,\, \lambda\in \mcal{P}(r),\, m_1(\lambda)=0\}$. Again, we take an orthogonal basis of the linear span of $\{s^r_\lambda(z);\, r \in \N_0,\, \lambda\in \mcal{P}(r),\, m_1(\lambda)=0\}$ consisting of weight vectors, so that each basis vector generates an
$\mfrak{sl}(2,\C)$-submodule isomorphic to $L_\mfrak{b}(\lambda \omega)$ for some $\lambda \in \C$ and all
these $\mfrak{sl}(2,\C)$-submodules are mutually orthogonal. Hence, we obtain
\begin{align}
  F \supset M\simeq  L_\mfrak{b}(0) \oplus M_\mfrak{b}(-2\omega) \oplus \bigoplus_{j=2}^\infty \bigoplus_{\substack{r_2,r_3,\dots \in \N_0 \\ 2r_2+3r_3+ \dots = j}} M_b(-2j\omega), \label{eq:submodule decomposition c=0}
\end{align}
where the direct summands are generated by $s^0_\emptyset(z)$, $s^1_\emptyset(z)$ and $\{s^r_\lambda(z);\, r\in \N_0,\, \lambda \in \mcal{P}(r),\, m_1(\lambda)=0,\, r+|\lambda|=j\}$ for $j\geq 2$.
\medskip

Let us recall that the character of the Verma module $M_\mfrak{b}(\lambda\omega)$ with highest weight $\lambda\omega \in \mfrak{h}^*$ for $\lambda \in \C$ is given by
\begin{align}
  \ch M_\mfrak{b}(\lambda \omega)= {e^{\lambda \omega} \over 1 - e^{-\alpha}} = {q^{-{1\over 2}\lambda} \over 1- q},
\end{align}
where $q=e^{-\alpha}=e^{-2\omega}$ is a formal variable.
\medskip

\theorem{\label{thm:module decomposition}
Let $\sigma_c \colon \mfrak{sl}(2,\C) \rarr \mfrak{gl}(F)$ be the representation of $\mfrak{sl}(2,\C)$ on $F$. Then we have
\begin{enumerate}
\item[1)] for $c \neq 0$
\begin{align}
  F \simeq  M_\mfrak{b}(-|c|^2\omega) \oplus \bigoplus_{j=2}^\infty \bigoplus_{\substack{r_2,r_3,\dots \in \N_0 \\ 2r_2+3r_3+ \dots = j}} M_\mfrak{b}(-(2j+|c|^2)\omega), \label{eq:module decomposition c not 0}
\end{align}
\item[2)] for $c=0$
\begin{align}
  F \simeq  L_\mfrak{b}(0) \oplus M_\mfrak{b}(-2\omega) \oplus \bigoplus_{j=2}^\infty \bigoplus_{\substack{r_2,r_3,\dots \in \N_0 \\ 2r_2+3r_3+ \dots = j}} M_\mfrak{b}(-2j\omega) \label{eq:module decomposition c=0}
\end{align}
\end{enumerate}
as $\mfrak{sl}(2,\C)$-modules. Recall the notation $L_\mfrak{b}(0)\simeq\,{\mathbb C}$
for the trivial simple $\mfrak{sl}(2,\C)$-module.}

\proof{Let us denote the right hand side of \eqref{eq:module decomposition c not 0} and \eqref{eq:module decomposition c=0} by $M$ depending on the value of a parameter $c \in \C$. Since $M \subset F$, as follows from \eqref{eq:submodule decomposition c not 0} and \eqref{eq:submodule decomposition c=0}, it is enough to compare the formal characters of $F$ and $M$. For the $\mfrak{sl}(2,\C)$-module $F$ we have
\begin{align*}
  \ch F &= \sum_{j=0}^\infty (\dim F^j)\,e^{-(2j+|c|^2)\omega} = q^{{1\over 2}|c|^2} \sum_{j=0}^\infty (\dim F^j)\,q^j = q^{{1\over 2}|c|^2} \sum_{j=0}^\infty \sum_{\substack{r_0,r_1,\dots \in \N_0 \\ r_0+2r_1+3r_2+\dots = j}} q^j \\
  &= q^{{1\over 2}|c|^2} \prod_{j=1}^\infty \sum_{k=0}^\infty (q^j)^k  = q^{{1\over 2}|c|^2} \prod_{j=1}^\infty {1 \over 1-q^j} = {q^{{1\over 2}|c|^2} \over \varphi(q)},
\end{align*}
where $\varphi(q)= \prod_{j=1}^\infty (1-q^j)$ and $F^j$ is the vector subspace of $F$ of graded homogeneous (with respect to the grading $\deg z_n =n+1$ for $n \in \N_0$) polynomials of homogeneity $j \in \N_0$. On the other hand, for the $\mfrak{sl}(2,\C)$-module $M$ we have
\begin{align*}
  \ch M &=
  \begin{cases}
    \ch M_\mfrak{b}(-|c|^2\omega) + \sum\limits_{j=2}^\infty \sum\limits_{\substack{r_2,r_3,\dots \in \N_0 \\ 2r_2+3r_3+ \dots = j}} \ch M_\mfrak{b}(-(2j+|c|^2)\omega)  & \text{for $c \neq 0$}, \\[5mm]
    \ch L_\mfrak{b}(0) + \ch M_\mfrak{b}(-2\omega) + \sum\limits_{j=2}^\infty \sum\limits_{\substack{r_2,r_3,\dots \in \N_0 \\ 2r_2+3r_3+ \dots = j}} \ch M_\mfrak{b}(-2j\omega) & \text{for $c=0$}.
  \end{cases}
\end{align*}
Using $\ch L_\mfrak{b}(0) + \ch M_\mfrak{b}(-2\omega) = \ch M_\mfrak{b}(0)$, we may write
\begin{align*}
  \ch M & = \ch M_\mfrak{b}(-|c|^2\omega) + \sum_{j=2}^\infty \sum_{\substack{r_2,r_3,\dots \in \N_0 \\ 2r_2+3r_3+ \dots = j}} \ch M_\mfrak{b}(-(2j+|c|^2)\omega) \\
  &= {q^{{1\over 2}|c|^2} \over 1-q} + \sum_{j=2}^\infty \sum_{\substack{r_2,r_3,\dots \in \N_0 \\ 2r_2+3r_3+ \dots = j}} {q^{j+{1\over 2}|c|^2} \over 1-q} = q^{{1\over 2}|c|^2} \bigg({1 \over 1-q} + \sum_{j=2}^\infty \sum_{\substack{r_2,r_3,\dots \in \N_0 \\ 2r_2+3r_3+ \dots = j}} {q^j \over 1-q}\bigg) \\
  &=  q^{{1\over 2}|c|^2} \bigg(\sum_{j=0}^\infty \sum_{\substack{r_2,r_3,\dots \in \N_0 \\ 2r_2+3r_3+ \dots = j}} {q^j \over 1-q}\bigg) = {q^{{1\over 2}|c|^2} \over 1-q}\, \prod_{j=2}^\infty \sum_{k=0}^\infty (q^j)^k \\
  &= {q^{{1\over 2}|c|^2} \over 1-q}\, \prod_{j=2}^\infty {1 \over 1 - q^j} = {q^{{1\over 2}|c|^2} \over \varphi(q)},
\end{align*}
which implies the statement.}

In fact, we have proved the following first main theorem, which gives an explicit description of the kernel of the linear differential operator $T_1^c$ on $F$ for all $c \in \C$.
\medskip

\theorem{\label{thm:kernel T_1^c}
The set
\begin{align}
  B_c =
  \begin{cases}
    \{t^c_\lambda(z);\, \lambda \in \mcal{P},\, m_1(\lambda)=0\} & \text{for $c\neq 0$}, \\
    \{s^r_\lambda(z);\, r \in \N_0,\, \lambda \in \mcal{P}(r),\, m_1(\lambda)=0\} & \text{for $c=0$}
  \end{cases}
\end{align}
forms a basis of $\ker T_1^c$ for all $c \in \C$.}

Now, we can formulate our second main theorem, which describes explicitly all
singular vectors in the Verma module $M_\ell(\delta,p)$ for the conformal Galilei algebra $\mfrak{cga}_\ell(1,\C)$ with $\ell \in \N$.
\medskip

\theorem{Let $\delta,p \in \C$ and $p \neq 0$. Then the vector space of singular vectors in the Verma module $M_\ell(\delta,p)$ for the conformal Galilei algebra $\mfrak{cga}_\ell(1,\C)$ with $\ell \in \N$ is a linear span of singular vectors
\begin{align}
  t^{{p\over \ell!}}_\lambda\!\bigg(\!-\!{P_{\ell-1} \over (\ell+1)!},-{P_{\ell-2} \over (\ell+2)!},\dots,-{P_0 \over (2\ell)!}\bigg)v_{\delta,p}
\end{align}
with character $(\delta+2|\lambda|)\widetilde{\omega}_D+p\widetilde{\omega}_{P_\ell}$,
where $\lambda \in \mcal{P}(\ell)$ with $m_1(\lambda)=0$. Here $v_{\delta,p}$ is the lowest
weight vector of $M_\ell(\delta,p)$ with lowest weight $\delta \omega_D + p \omega_{P_\ell}$,
and the polynomials $t^c_\lambda(z) \in F$ are defined by
\begin{align}
  t^c_\lambda(z)= \prod_{i \in \N} t^c_i(z)^{m_i(\lambda)}
\end{align}
with
\begin{align}
\begin{aligned}
  t^c_k(z_0,z_1,\dots,z_{k-1})&=(-1)^kc^k  u_k\Big(\!-\!{z_0 \over c},\dots,-{z_{k-1}\over c}\Big) \\ &= \sum_{\substack{r_1+2r_2+\dots+kr_k =k \\ r_1,r_2,\dots,r_k \in \N_0}} {k(r_1+\dots +r_k-1)! \over r_1!\dots r_k!} \prod_{i=1}^k (c^{i-1}z_{i-1})^{r_i}
\end{aligned}
\end{align}
for $c \in \C$.}

\proof{The statement is a straightforward consequence of Theorem \ref{thm:kernel T_1^c}, considerations after Lemma \ref{soloperz}, the transformation formula \eqref{eq:coord change} and the symmetrization map in Theorem \ref{thm:symmetrization}.}

The formulation of an analogous result for $p=0$ is much more delicate as follows from the subsequent example.
The explicit form of polynomials $s^r_\lambda(z)$ in \eqref{eq:s^r_lambda(z) example} implies
$\smash{s^3_{(2^3)}(z)}, \smash{s^3_{(3^2)}(z)} \in F_7$. However, we obtain
$\smash{4s^3_{(2^3)}(z)-3s^3_{(3^2)}(z)}  \in F_5$, which implies that
\begin{align*}
  s^3_{(2^3)}\!\bigg(\!-\!{P_{\ell-1} \over (\ell+1)!},-{P_{\ell-2} \over (\ell+2)!},\dots,-{P_0 \over (2\ell)!}\bigg)v_{\delta,0}, \qquad s^2_{(3^2)}\!\bigg(\!-\!{P_{\ell-1} \over (\ell+1)!},-{P_{\ell-2} \over (\ell+2)!},\dots,-{P_0 \over (2\ell)!}\bigg)v_{\delta,0}
\end{align*}
are singular vectors in $M_\ell(\delta,0)$ for $\ell \geq 7$ with a character $(\delta+18)\widetilde{\omega}_D$,
although
\begin{align*}
  (4s^3_{(2^3)}-3s^3_{(3^2)})\!\bigg(\!-\!{P_{\ell-1} \over (\ell+1)!},-{P_{\ell-2} \over (\ell+2)!},\dots,-{P_0 \over (2\ell)!}\bigg)v_{\delta,0}
\end{align*}
is a singular vector in $M_\ell(\delta,0)$ already for $\ell \geq 5$ with character $(\delta+18)\widetilde{\omega}_D$.
We may at least claim that
\begin{align}
  s^r_\lambda\!\bigg(\!-\!{P_{\ell-1} \over (\ell+1)!},-{P_{\ell-2} \over (\ell+2)!},\dots,-{P_0 \over (2\ell)!}\bigg)v_{\delta,0}
\end{align}
is a singular vector in $M_\ell(\delta,p)$ for $\ell \geq \ell_0$, where $\ell_0$ is the smallest natural number such that $s^r_\lambda(z) \in F_{\ell_0}$, with character $(\delta+2|\lambda|+2r)\widetilde{\omega}_D$ for $r \in \N_0$, $\lambda \in \mcal{P}(r)$ with $m_1(\lambda)=0$ and $v_{\delta,0}$ the lowest weight vector of $M_\ell(\delta,0)$ with lowest weight $\delta \omega_D$.
\medskip

We give explicit formulas for a few polynomials $s^r_\lambda(z)$ of small degree with $m_1(\lambda)=0$. We have
\begin{align}
  s^r_\emptyset(z)=z_0^r
\end{align}
for all $r \in \N_0$, and
\begin{align} \label{eq:s^r_lambda(z) example}
\begin{aligned}
  s^2_{(2^a)}(z)&= {1 \over 2^a} \sum_{\substack{i,j \in \N_0 \\ i+j=2a}} (-1)^i z_i z_j, \\
  s^3_{(2^a)}(z)&= {1 \over 2^a} \sum_{\substack{i,j \in \N_0 \\ i+j=2a}} (-1)^i z_0z_i z_j = z_0s^2_{(2^a)}, \\
  s^3_{(3^a)}(z)&= {1 \over 3^a} \sum_{\substack{i,j,k \in \N_0 \\ i+j+k = 3a}} \Big(\!\textstyle{-{1 \over 2}+ {\sqrt{3} \over 2}\,{\rm i}}\Big)^{i+2j} z_iz_jz_k
\end{aligned}
\end{align}
for all $a \in \N_0$, which easily follows from \eqref{eq:expansion} by a suitable choice of parameters $\alpha_k \in \C$ for $k=1,2\dots,r$, thus $\alpha_k=0$ for $k=1,2,\dots,r$ in the case $s^r_\emptyset(z)$, $\alpha_1=1$, $\alpha_2=-1$ for $\smash{s^2_{(2^a)}(z)}$, $\alpha_1=1$, $\alpha_2=-1$, $\alpha_3=0$ for $\smash{s^3_{(2^a)}(z)}$ and $\alpha_1=\omega$, $\alpha_2=\omega^2$, $\alpha_3=\omega^3$ for $\smash{s^3_{(3^a)}(z)}$, where $\omega=\smash{-{1 \over 2}+ {\sqrt{3} \over 2}\,{\rm i}}$ is the third root of unity.


\subsection{Singular vectors and differential equations of flag type}

The first order linear differential operator $T^c_1$ of the form \eqref{eq:operator T^c_1} is a representative of the class of partial
differential equations of flag type. A partial differential equation of flag type on a commutative $\C$-algebra $\mcal{A}=\C[x_1,x_2,\dots,x_n]$ is a linear partial differential equation of the form
\begin{align}
\big(D_1+f_1D_2+f_2D_3+\dots +f_{m-1}D_m\big)u=0 \label{eq:flag type}
\end{align}
for $u \in \mcal{A}$, where $D_1,D_2,\dots,D_m$ for $m \in \N$ are commuting locally nilpotent algebraic differential operators on the $\C$-algebra $\mcal{A}$ and $f_1,f_2,\dots,f_{m-1} \in \mcal{A}$ are polynomials satisfying
\begin{align}
  D_k(f_j)=0\quad \text{for} \quad k > j.
\end{align}
The results of \cite{Xu2008} and \cite{Xu2013} gives a parametrization of solution spaces of partial differential equations of flag
type \eqref{eq:flag type} in the following form. Let $\mcal{B}$ be a $\C$-subalgebra of $\mcal{A}$, so that $\mcal{A}$ is a free $\mcal{B}$-module generated by a filtered vector subspace $\mcal{V} = \bigcup_{r=0}^\infty \mcal{V}_r$, $\mcal{V}_r \subset \mcal{V}_{r+1}$ for $r \in \N_0$, of $\mcal{A}$. Furthermore, let $T_1$ be a linear differential operator on $\mcal{A}$ with a right inverse $T_1^-$ such that
\begin{align}
  T_1(\mcal{B}) \subset \mcal{B}, \qquad T_1^-(\mcal{B}) \subset \mcal{B}, \qquad T_1(bv)=T_1(b)v, \qquad T_1^-(bv)=T_1^-(b)v
\end{align}
for all $b \in \mcal{B}$ and $v \in \mcal{V}$, and let $T_2$ be a linear differential operator on $\mcal{A}$ such that
\begin{align}
  T_2(\mcal{V}_r) \subset \mcal{B}\mcal{V}_{r-1}, \qquad T_2(ba)=bT_2(a)
\end{align}
for all $a \in \mcal{A}$, $b \in \mcal{B}$ and $r \in \N_0$ with the notation $\mcal{V}_{-1}=\{0\}$. Then the vector subspace of all $u\in \mcal{A}$
fulfilling $(T_1+T_2)(u)=0$ is given by the complex linear span
\begin{align}
\Big\langle \Big\{{\textstyle \sum\limits_{j=0}^\infty} (-T_1^- T_2)^j(bv);\, v \in \mcal{V},\, b \in \mcal{B},\, T_1(b)=0\Big\}\Big\rangle. \label{eq:solution flag type}
\end{align}
In fact, there are only finitely many contributions in the last summation formula.
\medskip

In the case of the linear differential operator $T_1^c$, $c \neq 0$, which acts on the $\C$-algebra $\mcal{A}=F_\ell=\C[z_0,z_1,\dots,z_{\ell-1}]$ as
\begin{align}
  - c\partial_{z_0}+\sum_{n=0}^{\ell-2} z_n \partial_{z_{n+1}},
\end{align}
we take $\mcal{B}=\C[z_0]$ and $\mcal{V}=\C[z_1,z_2,\dots,z_{\ell-1}]$ with $\mcal{V}_r = \mcal{V} \cap \bigoplus_{k=0}^r \mcal{A}_k$, where $\mcal{A}_k \subset \mcal{A}$ is the subspace of graded homogeneous polynomials of homogeneity $k$ with respect to the grading $\deg z_n = n+1$ for $n=0,1,\dots,\ell-1$. Furthermore, if we set
\begin{align}
  T_1 = -c\partial_{z_0}\qquad \text{and} \qquad T_2= \sum\limits_{n=0}^{\ell-2} z_n \partial_{z_{n+1}},
\end{align}
then its right inverse $T_1^-$ is given by
\begin{align}
  T_1^-(f(z_0,z_1,\dots,z_{\ell-1}))= -{1 \over c} \int_0^{z_0}f(z'_0,z_1,\dots,z_{\ell-1})\, {\rm d}z'_0,
\end{align}
where $\int_0^{z_0}\!f(z'_0,z_1,\dots,z_{\ell-1})\,{\rm d}z'_0$ denotes the definite integration in the variable $z'_0$. In conclusion, we obtain
\begin{align}
 \ker\!\bigg(\!-\!c\partial_{z_0} + \sum_{n=0}^{\ell-2} z_n \partial_{z_n}\!\bigg) =  \ker(T_1+T_2)    = \bigg\{\sum_{j=0}^\infty (-T_1^-T_2)^j(v);\, v \in \mcal{V}\bigg\}
\end{align}
for $c \neq 0$.


\begin{appendices}

\section{Symmetric polynomials and functions}
\label{app:symetric polynomials and fucntions}

For the reader's convenience, we summarize several basic facts concerning symmetric polynomials and symmetric functions used throughout the article, see \cite{Macdonald1995} for a detailed information.
\medskip

Let $S_n$ be the symmetric group on the set $\{1,2,\dots,n\}$ and $\C^n$ its defining representation. Then we have the induced representation of $S_n$ on the vector space $\C[\C^n]$ of polynomials on $\C^n$ defined by
\begin{align}
  (\sigma \cdot f)(v)= f(\sigma^{-1}v),
\end{align}
where $\sigma \in S_n$, $f \in \C[\C^n]$ and $v \in \C^n$. The subset $\C[\C^n]^{S_n}$ of all invariant polynomials is a $\C$-subalgebra of $\C[\C^n]$, and is called the ring of symmetric polynomials.

Let us denote by $x=(x_1,x_2,\dots,x_n)$ the canonical linear coordinate functions on $\C^n$. Then we have an isomorphism $\C[\C^n] \simeq \C[x]$, and the induced representation of $S_n$ on $\C[x]$ is given by
\begin{align}
  \sigma \cdot f(x_1,x_2,\dots,x_n) = f(x_{\sigma(1)},x_{\sigma(2)},\dots,x_{\sigma(n)})
\end{align}
for $\sigma \in S_n$ and $f(x_1,x_2,\dots,x_n) \in \C[x]$. We denote by $\Lambda_n(x)$ the ring of $S_n$-invariant polynomials in the variables $x_1,x_2,\dots,x_n$. Moreover, the ring $\Lambda_n(x)$ is a graded $\C$-algebra, we have
\begin{align}
  \Lambda_n(x) = \bigoplus_{k=0}^\infty \Lambda_n^k(x),
\end{align}
where $\Lambda_n^k(x)$ is the vector space of homogeneous symmetric polynomials of degree $k$.
\medskip

The most important bases of the ring of symmetric polynomials are labeled by partitions of non-negative integers $\N_0$. A partition of $n \in \N_0$ is a non-increasing sequence $\lambda=(\lambda_1,\lambda_2,\dots,\lambda_m)$ of positive integers called parts whose sum is $n$, i.e.\ $\lambda_1 \geq \lambda_2 \geq \dots \geq \lambda_m >0$ and $n=\sum_{i=1}^m \lambda_i$. The only partition of $0$ is the empty partition $\emptyset$. The length $\ell(\lambda)=m$ of $\lambda$ is the length of the sequence and the size $|\lambda|=\sum_{i=1}^m\lambda_i$ of $\lambda$ is the sum over all elements in the sequence. It is often convenient to define $\lambda_i=0$ for $i >\ell(\lambda)$. The set of all partitions of $n$ is denoted by $\mcal{P}_n$, the set of all partitions by $\mcal{P}$, the set of all partitions of $n$ whose parts are at most $r$ by $\mcal{P}_n(r)$, and the set of all partitions whose parts are at most $r$ by $\mcal{P}(r)$. Furthermore, for a partition $\lambda$ we define
\begin{align}
  \veps_\lambda = (-1)^{|\lambda|-\ell(\lambda)},\qquad  \quad z_\lambda = \prod_{i \in \N} i^{m_i(\lambda)}m_i(\lambda)!,\qquad  \quad u_\lambda= {\ell(\lambda)! \over \prod_{i\in \N}m_i(\lambda)!},
\end{align}
where $m_i(\lambda)$ indicates the number of times the integer $i$ occurs in the partition $\lambda$, it is called the multiplicity of $i$ in $\lambda$. We shall also use a notation $(1^{m_1(\lambda)},2^{m_2(\lambda)},\dots,r^{m_r(\lambda)},\dots)$ for a partition $\lambda$.


\subsection{Symmetric polynomials}

Here we recall several important classes of symmetric polynomials which are useful in the construction of generating sets of $\Lambda_n(x)$.

\begin{enumerate}
{\bf \item[1)] Monomial symmetric polynomials.} For a partition $\lambda$ satisfying $\ell(\lambda) \leq n$, we define the symmetric polynomial
\begin{align}
  m_\lambda(x_1,\dots,x_n)= \sum_{\{i_1,i_2,\dots,i_{\ell(\lambda)}\} \subset \{1,2\dots,n\}} x_{i_1}^{\lambda_1}x_{i_2}^{\lambda_2}\dots x_{i_{\ell(\lambda)}}^{\lambda_{\ell(\lambda)}}
\end{align}
called monomial symmetric polynomial.

{\bf \item[2)] Elementary symmetric polynomials.} The symmetric polynomials
\begin{align}
  e_k(x_1,\dots,x_n) =
  \begin{cases}
   1 & \text{for $k=0$}, \\
   \sum\limits_{1\leq i_1 < i_2 < \dots < i_k \leq n} x_{i_1}x_{i_2}\dots x_{i_k} &  \text{for $k \geq 1$}
  \end{cases}
\end{align}
are called elementary symmetric polynomials. For a partition $\lambda$, we define
\begin{align}
  e_\lambda(x_1,\dots,x_n)= \prod_{i=1}^{\ell(\lambda)} e_{\lambda_i}(x_1,\dots,x_n).
\end{align}

{\bf \item[3)] Power sum symmetric polynomials.} The symmetric polynomials
\begin{align}
  p_k(x_1,\dots,x_n) = \sum_{i=1}^n x_i^k \qquad \text{for $k \geq 1$}
\end{align}
are called power sum symmetric polynomials. For a partition $\lambda$, we define
\begin{align}
  p_\lambda(x_1,\dots,x_n)= \prod_{i=1}^{\ell(\lambda)} p_{\lambda_i}(x_1,\dots,x_n).
\end{align}

{\bf \item[4)] Complete homogeneous symmetric polynomials.} The symmetric polynomials
\begin{align}
  h_k(x_1,\dots,x_n)= \begin{cases}
  1 & \text{for $k=0$}, \\
  \sum\limits_{1 \leq i_1 \leq i_2 \leq \dots \leq i_k \leq n} x_{i_1}x_{i_2}\dots x_{i_k} & \text{for $k \geq 1$}
  \end{cases}
\end{align}
are called complete homogeneous symmetric polynomials. For a partition $\lambda$, we define
\begin{align}
  h_\lambda(x_1,\dots,x_n)= \prod_{i=1}^{\ell(\lambda)} h_{\lambda_i}(x_1,\dots,x_n).
\end{align}
\end{enumerate}

Moreover, there is an isomorphism of $\C$-algebras
  \begin{align}
    \Lambda_n(x) \simeq \C[e_1(x_1,\dots,x_n),\dots,e_n(x_1,\dots,x_n)] \simeq \C[p_1(x_1,\dots,x_n),\dots,p_n(x_1,\dots,x_n)]
  \end{align}
for all $n\in \N$.


\subsection{Newton's identities}

We review the Newton's identities, also known as the Newton-Girard formulae, which give relations between power sum and elementary symmetric polynomials.

It is more convenient to work with symmetric polynomials in infinitely than finitely many variables. For any pair $m,n\in \N$ with $m \geq n$, we consider the $\C$-algebra homomorphism
\begin{align}
  \rho_{m,n} \colon \Lambda_m(x) \rarr \Lambda_n(x)
\end{align}
defined by
\begin{align}
  \rho_{m,n}(x_i) = \begin{cases}
    x_i & \text{for $i=1,\dots,n$}, \\
    0   & \text{for $i=n+1,\dots,m$}.
  \end{cases}
\end{align}
Then the ring of symmetric functions $\Lambda(x)$ in infinitely many variables $x_n$, $n \in \N$, is defined by the inverse limit
\begin{align}
  \Lambda(x) = \lim_{\longleftarrow n} \Lambda_n(x)
\end{align}
together with the projections $\rho_n \colon \Lambda(x) \rarr \Lambda_n(x)$ given by
\begin{align}
  \rho_n(x_i) = \begin{cases}
    x_i & \text{for $i=1,\dots,n$}, \\
    0   & \text{for $i > n$}.
  \end{cases}
\end{align}
Let us denote by $e_n(x)$ and $p_n(x)$ symmetric functions in $\Lambda(x)$ satisfying
\begin{align}
  \rho_k(e_n(x))=e_n(x_1,x_2,\dots,x_k) \quad \text{and} \quad \rho_k(p_n(x))=p_n(x_1,x_2,\dots,x_k)
\end{align}
for all $k \in \N$, respectively. As above, for a partition $\lambda$ we define
\begin{align}
  e_\lambda(x) = \prod_{i=1}^{\ell(\lambda)} e_{\lambda_i}(x) \quad \text{and} \quad p_\lambda(x) = \prod_{i=1}^{\ell(\lambda)} p_{\lambda_i}(x).
\end{align}
For $n \in \N$, the transformation between elementary and power sum symmetric polynomials is
\begin{align}
  e_n(x)= \sum_{\lambda \in \mcal{P}_n} \veps_\lambda {p_\lambda(x) \over z_\lambda}=  \sum_{\substack{r_1+2r_2+\dots+nr_n =n \\ r_1,r_2,\dots, r_n \in \N_0}} (-1)^n \prod_{i=1}^n {(-p_i(x))^{r_i} \over r_i! i^{r_i}} \label{eq:expression e in p}
\end{align}
and
\begin{align}
  p_n(x)= n\sum_{\lambda \in \mcal{P}_n} {\veps_\lambda u_\lambda \over \ell(\lambda)}\, e_\lambda(x) = \sum_{\substack{r_1+2r_2+\dots+nr_n =n \\ r_1,r_2,\dots,r_n \in \N_0}} (-1)^n {n(r_1+\dots +r_n-1)! \over r_1!\dots r_n!} \prod_{i=1}^n (-e_i(x))^{r_i}. \label{eq:expression p in e}
\end{align}
Moreover, if we define polynomials
\begin{align}
  u_n(y_1,y_2,\dots,y_n)=\sum_{\substack{r_1+2r_2+\dots+nr_n =n \\ r_1,r_2,\dots,r_n \in \N_0}} (-1)^n {n(r_1+\dots +r_n-1)! \over r_1!\dots r_n!} \prod_{i=1}^n (-y_i)^{r_i}  \label{eq:u_n polynomial}
\end{align}
for $n \in \N$, then we may write
\begin{align}
  p_n(x)=u_n(e_1(x),e_2(x),\dots,e_n(x))  \label{eq:relation p and e}
\end{align}
for all $n \in \N$. Since we have $e_r(x_1,\dots,x_n)=0$ for $r>n$, we obtain from \eqref{eq:expression p in e} for $k>n$ the formula
\begin{align}
\begin{aligned}
  p_k(x_1,\dots,x_n)&= k\sum_{\lambda \in \mcal{P}_k\!(n)} {\veps_\lambda u_\lambda \over \ell(\lambda)}\, e_\lambda(x_1,\dots,x_n) \\
  &=\sum_{\substack{r_1+2r_2+\dots+nr_n =k \\ r_1,r_2,\dots, r_n \in \N_0}} (-1)^k {k(r_1+r_2\dots +r_n-1)! \over r_1!r_2!\dots r_n!} \prod_{i=1}^n (-e_i(x_1,\dots,x_n))^{r_i}.
\end{aligned}  \label{eq:expression p in e polynomials}
\end{align}
The equality \eqref{eq:expression e in p} can be used to produce $p_k(x_1,\dots,x_n)\in \C[p_1(x_1,\dots,x_n),\dots,p_n(x_1,\dots,x_n)]$  as a polynomial in $p_1(x_1,\dots,x_n),\dots,p_n(x_1,\dots,x_n)$ for $k>n$.
\medskip

The expression of $p_k(x_1,\dots,x_n)$ for all $k \in \N$ in the quotient ring
\begin{align}
\Lambda_n(x)/(p_1(x_1,\dots,x_n)) \simeq \C[p_2(x_1,\dots,x_n),\dots,p_n(x_1,\dots,x_n)]
\end{align}
easily follows from \eqref{eq:expression p in e polynomials}. If we restrict to small values of $n$, we obtain the following identities.
\begin{enumerate}
  \item[1)] If $n=1$, then we have
  \begin{align}
    p_k(x_1) = 0
  \end{align}
  for $k \geq 1$ in $\Lambda_1(x)/(p_1(x_1))$.
  \item[2)] If $n=2$, then we have
  \begin{align}
    p_k(x_1,x_2) = \begin{cases}
      2^{1-{k \over 2}} (p_2(x_1,x_2))^{{k \over 2}}  & \text{for $k\geq 1$ even}, \\
      0  & \text{for $k\geq 1$ odd}
    \end{cases}
  \end{align}
  in $\Lambda_2(x)/(p_1(x_1,x_2))$.
  \item[3)] If $n=3$, then we have
  \begin{align}
    p_k(x_1,x_2,x_3)= \sum_{\substack{2r_2+3r_3=k \\ r_2,r_3 \in \N_0}} {2r_2+3r_3 \over r_2+r_3} \binom{r_2+r_3}{r_2} {(p_2(x_1,x_2,x_3))^{r_2} (p_3(x_1,x_2,x_3))^{r_3} \over 2^{r_2} 3^{r_3}}
  \end{align}
  for $k \geq 1$ in $\Lambda_3(x)/(p_1(x_1,x_2,x_3))$.
\end{enumerate}


\subsection{The Cauchy kernel}

Let $x=(x_1,x_2,\dots)$ and $y=(y_1,y_2,\dots,)$ be two infinite sequences of independent variables. We shall review the series expansions for the infinite product
\begin{align}
  \Omega^* = \prod_{i,j \in \N}(1+x_i y_j),
\end{align}
called the dual Cauchy kernel.
\medskip

\proposition{We have
\begin{align}
  \prod_{i,j \in \N} (1+x_iy_j)= \sum_{\lambda\in \mcal{P}} e_\lambda(x) m_\lambda(y) = \exp\!\bigg(\sum_{n=1}^\infty {(-1)^{n+1} \over n}\, p_n(x)p_n(y)\!\bigg) = \sum_{\lambda \in \mcal{P}} \veps_\lambda {p_\lambda(x) p_\lambda(y) \over z_\lambda}. \label{eq:Cauchy kernel}
\end{align}}

\proof{For the first identity in \eqref{eq:Cauchy kernel}, we may write
\begin{align*}
  \prod_{i,j\in \N}(1+x_iy_j)&=\prod_{j=1}^\infty \bigg(\prod_{i=1}^\infty (1+x_iy_j) \!\bigg) =  \prod_{j=1}^\infty \sum_{k=0}^\infty e_k(x)y_j^k = \sum_{\lambda \in \mcal{P}} e_\lambda(x) m_\lambda(y).
\end{align*}
For the second identity in \eqref{eq:Cauchy kernel}, let us recall the power series expansions of
\begin{align*}
  \log(1+z) = \sum_{n=1}^\infty (-1)^{n+1} {z^n \over n} \qquad \text{and} \qquad \exp(z)= \sum_{n=0}^\infty {z^n \over n!}.
\end{align*}
Therefore, we have
\begin{align*}
  \log \prod_{i,j \in \N} (1+x_i y_j) &= \sum_{i,j \in \N} \log(1+x_iy_j) = \sum_{i,j\in \N} \sum_{n=1}^\infty (-1)^{n+1} {(x_iy_j)^n \over n} \\
  & = \sum_{n=1}^\infty {(-1)^{n+1} \over n} \bigg(\sum_{i=1}^\infty x_i^n \bigg)\!\bigg(\sum_{j=1}^\infty y_j^n \bigg)= \sum_{n=1}^\infty {(-1)^{n+1} \over n}\, p_n(x)p_n(y),
\end{align*}
which implies
\begin{align*}
  \prod_{i,j\in \N} (1+x_i y_j) &= \exp\!\bigg(\sum_{n=1}^\infty {(-1)^{n+1} \over n}\, p_n(x)p_n(y)\! \bigg) = \sum_{k=0}^\infty {1 \over k!} \bigg(\sum_{n=1}^\infty {(-1)^{n+1} \over n}\, p_n(x)p_n(y)\! \bigg)^k \\
  &= \sum_{k=0}^\infty {1 \over k!} \bigg( \sum_{r_1,r_2,\dots \in \N_0}\! \binom{k}{r_1,r_2,\dots}\!\bigg({p_1(x)p_1(y) \over 1} \bigg)^{\!r_1}\!\bigg(\!-{p_2(x)p_2(y) \over 2}\bigg)^{\!r_2}\! \dots \bigg)\\
  & = \sum_{\lambda \in \mcal{P}} \veps_\lambda {p_\lambda(x) p_\lambda(y) \over z_\lambda},
\end{align*}
where $\binom{k}{r_1,r_2,\dots}= {k! \over r_1!r_2!\dots}$ for $r_1+r_2+\dots = k$ denotes the multinomial coefficient. This completes the proof.}

\end{appendices}

\vspace{-2mm}


\section*{Conclusion and outlook}

As we already noticed, conformal Galilei algebras $\mfrak{cga}_1(d,\R)$ are given by the Inönu-Wigner contraction of conformal Lie algebras $\mfrak{so}(d+1,2,\R)$. Homomorphisms of generalized Verma modules for conformal Lie algebras were already classified, see e.g.\ \cite{koss2015}, so it would be desirable to understand the process of contraction directly on the level of homomorphisms of (generalized) Verma modules (or singular vectors).

Secondly, the singular vectors constructed in the present article are in bijective correspondence with differential operators equivariant for a conformal Galilei algebra and acting on smooth sections of line bundles supported on the homogeneous space $G/B$. The explicit form of singular vectors allows to deduce the formulas for equivariant differential operators, acting in the non-compact picture of induced representations.

We also remark that the representation of $\mfrak{sl}(2,\C)$ on the vector space $F$, defined in Section \ref{sec:vermmodsingsec}, can be extended to the so called Fock representation of the Virasoro algebra, cf.\ \cite{Wakimoto-Yamada1986}.


\section*{Acknowledgement}

L.\,Křižka and P.\,Somberg acknowledge the financial support from the grant GA\,CR P201/12/G028.



\providecommand{\bysame}{\leavevmode\hbox to3em{\hrulefill}\thinspace}
\providecommand{\MR}{\relax\ifhmode\unskip\space\fi MR }
\providecommand{\MRhref}[2]{%
  \href{http://www.ams.org/mathscinet-getitem?mr=#1}{#2}
}
\providecommand{\href}[2]{#2}

\end{document}